\documentclass[11pt,reqno]{article}

\usepackage[text={153mm,235mm},centering]{geometry}                
\usepackage{booktabs}
\usepackage{tabularx}
\usepackage{graphicx}
\usepackage{amssymb,bm}
\usepackage{epstopdf}
\usepackage{array}
\usepackage{amsmath,amsthm}
\usepackage{amssymb} 
\usepackage[mathscr]{eucal}

\usepackage{comment}
\usepackage{color}
\usepackage[all]{xy}
\usepackage{mathrsfs}
\usepackage{authblk}
\usepackage{hyperref}

\newtheorem{theorem}{Theorem}[section]
\newtheorem{lemma}[theorem]{Lemma}

\newtheorem{corollary}[theorem]{Corollary}
\newtheorem{proposition}[theorem]{Proposition}

\newtheorem{definition-lemma}[theorem]{Definition-Lemma}
\newtheorem{definition-theorem}[theorem]{Definition-Theorem}

\theoremstyle{definition}
\newtheorem{example}[theorem]{Example}
\newtheorem{definition}[theorem]{Definition}

\newtheorem{remark}[theorem]{Remark}

\newtheorem*{ack}{Acknowledgements}

\newtheorem{convention}[theorem]{Convention}

\newcommand{\ac}{\textup{!`}}

\usepackage{marginnote}
\let\oldmarginpar\marginpar
\renewcommand\marginpar[1]{\-\oldmarginpar[\raggedleft\footnotesize
\color{red} #1]
{\raggedright\footnotesize\color{red} #1}} 
\allowdisplaybreaks

\begin{document}

\title{Differential calculus on bigraded spaces and Koszul duality}

\author[1,2]{Ruobing Chen\thanks{Email: astrocola@163.com}}

\author[3]{Sirui Yu\thanks{Email:
banaenoptera@163.com}}

\renewcommand\Affilfont{\small}

\affil[1]{School of Mathematics, Sichuan University, 
Chengdu 610064, P.R. China}

\affil[2]{Department of Mathematics, New Uzbekistan University, 
Tashkent 100001, Uzbekistan}

\affil[3]{School of Science,
Southwest Petroleum University,
Chengdu 610500, P.R. China}

\date{}

\maketitle
\begin{abstract}
Koszul duality for superalgebras involves both the intrinsic super
parity and the degree shift introduced by suspension. We study this
duality while keeping the two gradings separate, using a
$\mathbb Z\times\mathbb Z$-graded refinement of the commutativity
convention of Deligne and Morgan. In this setting, the Koszul
suspension preserves the parity degree and raises the independent
weight by one. We realize the polynomial superalgebra
$\mathcal O(k^{m|n})$ and its Koszul dual as function algebras on
bigraded spaces and construct an isomorphism between their
differential calculi of polyvector fields and mixed differential
forms. For admissible quadratic Poisson structures, whose bivectors
have bidegree $(0,2)$, this isomorphism descends to the differential
calculi on Poisson cohomology and mixed Poisson homology. With
respect to the corresponding mixed volume forms, Koszul duality
preserves unimodularity and, in the unimodular case, identifies
the induced Batalin--Vilkovisky algebra structures on Poisson
cohomology.
\end{abstract}

\setcounter{tocdepth}{3} \tableofcontents


\section{Introduction}\label{sect:intro}
Let $k^{m|n}$ be the usual superspace over a field $k$ of characteristic
zero. Its algebra of algebraic functions is the supercommutative
polynomial algebra
\begin{equation}\label{funonRmn}
\mathcal O(k^{m|n})
=
\mathrm{Sym}_k\{x_1,\ldots,x_m,y_1,\ldots,y_n\},
\end{equation}
where $|x_i|=0$ and $|y_j|=1$ for $i=1,\ldots,m$ and
$j=1,\ldots,n$. In the two extreme cases,
$\mathcal O(k^{m|0})$ is the ordinary polynomial algebra and
$\mathcal O(k^{0|m})$ is the exterior algebra; after identifying the
generating vector spaces with the appropriate duals, these two
algebras are Koszul dual to each other.
In particular, taking $m=n$, we get that $k^{m|m}$ is Koszul self-dual.
For an account of Koszul duality in quantum field theory, including its
role in coupling quantum field theories to topological line defects, see
Paquette and Williams~\cite{paquette2022}.

In practice, however, if we want to study the homological properties
of superspaces like $k^{m|n}$, some grading issue arises.
In \cite{deligne1999notes},
Deligne and Morgan proposed to equip variables
in a superspace with a bi-grading.
Their treatment
of K\"ahler differentials illustrates how the two gradings can be
handled separately: one may first define the universal derivation
as a parity-preserving map and then assign degrees to differential
forms. We use a related procedure in Section~\ref{sect:diffcalcul}:
the universal derivation on a bigraded algebra preserves both
degrees, while the differential on forms is obtained by a shift
of weight and has bidegree $(0,-1)$.
These examples motivate the
use of two separate gradings in the differential calculus and
Koszul duality studied here.

From the viewpoint of derived constructions in supergeometry,
keeping the integer grading separate allows the shifts occurring
in complexes to be recorded independently of the original super
parity. The graded super setting considered here also provides a
basic example for our subsequent work on two-brane formality and
deformation quantization with commutation factors. In that
setting, the internal commutation rule is generalized while the
cohomological grading remains separate.

We realize $\mathcal O(k^{m|n})$ and its Koszul dual as function
algebras on two $\mathbb Z\times\mathbb Z$-graded superspaces and
study their differential calculus and Poisson structures.
Throughout the paper, the two degrees are called the
\textit{parity} and the \textit{weight}. The first grading
refines the usual super parity, while the second is an independent
integer grading that participates in the commutativity rule and
records the shifts used in our constructions. Its relation to the
homological and cohomological indices of the complexes considered
below is specified in their definitions.
Reducing the parity degree modulo $2$ gives the
$\mathbb Z_2\times\mathbb Z$-grading and commutativity convention
of Deligne and Morgan. We retain the integer refinement in
order to record the finer homogeneity conditions used below.

\subsection{Koszul duality of bigraded
 spaces}

We first
identify $k^{m|n}$
as a $\mathbb{Z}\times\mathbb{Z}$-graded space. Its function algebra remains \eqref{funonRmn}, but now
$\mathrm{deg}(x_i)=(0,0)$
and $\mathrm{deg}(y_j)=(1,0)$.
We call the first component of the degree the \textit{parity degree} and the second the \textit{weight}; the signs depend on their reductions modulo $2$.
Following the sign convention
proposed by Deligne-Morgan in
\cite{deligne1999notes},
the commutativity
of the variables is given by
\begin{equation}\label{signconv}
x  y=(-1)^{p(x)p(y)+w(x)w(y)}y x,
\end{equation}
where $p(-)$ and $w(-)$ denote the parity degree and weight, respectively.

The suspension \(s\) used in quadratic duality has bidegree \((0,1)\): it preserves parity and raises the weight by one.
The Koszul dual algebra
of \eqref{funonRmn}
is given by
\begin{equation}\label{funonRmnKoszul}
\mathcal{O}(k^{n\wedge m}):=\mathrm{Sym}_{k}
\{ \xi_1,
\cdots, \xi_m,\eta_1,\cdots,
\eta_n\},
\end{equation}
where \(\xi_i=sx_i^*\) and \(\eta_j=sy_j^*\), with \(\deg(\xi_i)=(0,1)\) and \(\deg(\eta_j)=(-1,1)\), for \(i=1,\ldots,m\) and \(j=1,\ldots,n\).
The commutativity signs again follow from \eqref{signconv}.
In particular, the $\xi_i$
anticommute with one another, the $\eta_j$ commute with one
another, and $\xi_i\eta_j=-\eta_j\xi_i$.

The relation with the convention using parity reversal is given
by the change of convention explained in the appendix to
Section~1 of \cite{deligne1999notes}. Replacing parity by $p+w$
modulo $2$, together with the product change
$a*b=(-1)^{w(a)p(b)}ab$, gives an ordinary supercommutative algebra
in which the $\eta_j$ are even and the $\xi_i$ are odd. Thus this
passage involves both a regrading and the corresponding change of
multiplication signs. We retain the bigrading and the
multiplication governed by \eqref{signconv} throughout this paper.

To distinguish the spaces carrying these bigraded function
algebras from ordinary superspaces, we denote the spaces
underlying \eqref{funonRmn} and \eqref{funonRmnKoszul} by
$k^{m||n}$ and $k^{n\wedge m}$, respectively. We obtain the
following theorem.

\begin{theorem}[Theorem
\ref{newmainthm1}]\label{mainthm1}
With the above notations,
\eqref{funonRmn} and
\eqref{funonRmnKoszul} are Koszul
dual to each other
(in the sense of Definition
\ref{DefKoszul} below);
in other words,
$k^{m||n}$ and
$k^{n\wedge m}$
are Koszul dual spaces.
\end{theorem}

\subsection{The differential calculus
structure}

Now we study the differential calculus carried by polyvector fields and
mixed differential forms. On a smooth manifold, polyvector fields form a Gerstenhaber algebra
under the wedge product and the Schouten--Nijenhuis bracket
\cite{Koszul1985crochet}, whereas differential forms carry the
corresponding contraction and Lie derivative actions. The Cartan
identities express the compatibility of these operations with the de
Rham differential. Following Tamarkin and Tsygan
\cite{tamarkin2005ring}, we refer to this collection of compatible
operations as a differential calculus.

We establish a bigraded analogue of this structure and then compare the
differential calculi associated with the Koszul dual algebras
$A=\mathcal O(k^{m||n})$,
$A^!=\mathcal O(k^{n\wedge m})$.
On the generators of the polyvector field and mixed form algebras, the
relevant bidegree preserving assignments are
\begin{equation}\label{eq:intro-calculus-correspondence}
\begin{aligned}
\phi_{\mathfrak X}:\quad&
x_i\longmapsto\partial_{\xi_i},
& y_j\longmapsto\partial_{\eta_j},
&\quad \partial_{x_i}\longmapsto-\xi_i,
& \partial_{y_j}\longmapsto\eta_j,\\
\phi_{\Omega}:\quad&
x_i\longmapsto\partial_{\xi_i},
& \partial_{y_j}\longmapsto\eta_j,
&\quad dx_i\longmapsto\xi_i^*,
& y_j^*\longmapsto d\eta_j,
\end{aligned}
\end{equation}
where \(1\leq i\leq m\) and \(1\leq j\leq n\).  These assignments are
induced by an isomorphism between $\widetilde\Omega^0(A)$ and
$\widetilde\Omega^0(A^!)$;
their functorial compatibility is proved in
Proposition~\ref{prop:fourier-compatibility}.  The result is summarized in
the following theorem.

\begin{theorem}[Theorems~\ref{DC-on-bigraded}
and~\ref{thm:dc-koszul-duality}]
(1) Let \(\Sigma\) be 
a finite-dimensional bigraded vector space. Then
\[
\bigl(
\mathfrak X^{\bullet,\bullet}(\mathcal O(\Sigma)),
\widetilde\Omega^{\bullet,\bullet}(\mathcal O(\Sigma)),
\wedge,\widetilde\iota,[-,-],\widetilde d
\bigr)
\]
forms a differential calculus.

(2) For the Koszul dual pair 
\(A\) and \(A^!\) above, the assignments in
\eqref{eq:intro-calculus-correspondence} extend to maps
\(\phi_{\mathfrak X}\) and \(\phi_{\Omega}\) that define an isomorphism
of differential calculi
\[
\phi=(\phi_{\mathfrak X},\phi_{\Omega}):
\bigl(
\mathfrak X^{\bullet,\bullet}(A),
\widetilde\Omega^{\bullet,\bullet}(A)
\bigr)
\xrightarrow{\;\sim\;}
\bigl(
\mathfrak X^{\bullet,\bullet}(A^!),
\widetilde\Omega^{\bullet,\bullet}(A^!)
\bigr).
\]
\end{theorem}

\subsection{Koszul duality of
Poisson structures}
We next study Poisson structures on $k^{m||n}$ and
$k^{n\wedge m}$. These structures are encoded by bivectors $\pi$
satisfying $[\pi,\pi]_{\mathrm{SN}}=0$, with graded commutativity and all
Schouten--Nijenhuis signs governed by \eqref{signconv}. The quadratic
Poisson structures entering our comparison are homogeneous of bidegree
$(0,2)$.

On $\mathcal O(k^{m||n})$, write the generators
$\{x_1,\ldots,x_m,y_1,\ldots,y_n\}$ collectively as
$\{z_1,\ldots,z_{m+n}\}$.  Similarly, on
$\mathcal O(k^{n\wedge m})$, write
$\{\xi_1,\ldots,\xi_m,\eta_1,\ldots,\eta_n\}$ collectively as
$\{\zeta_1,\ldots,\zeta_{m+n}\}$.

A mixed bivector on $k^{m||n}$ is quadratic if it has the form
\begin{equation}\label{quadraticPoisson}
\pi=\sum_{i,j,p,q}c_{ij}^{pq}z_iz_j\partial_p\wedge\partial_q.
\end{equation}
 Under the
correspondence \eqref{eq:intro-calculus-correspondence}, $\pi$ is sent to
a quadratic bivector $\pi^!$ on $k^{n\wedge m}$, called its
\emph{Koszul dual}.  Proposition~\ref{quadratic poisson Koszul} proves that a quadratic bivector $\pi$ is Poisson if and only if $\pi^!$ is
Poisson.

For a smooth Poisson manifold, the Poisson cochain complex is
constructed from polyvector fields; see
\cite{lichnerowicz1977varietes}. The Poisson chain complex is
constructed from differential forms; see
\cite[Chapter 3]{laurent2012poisson}.  An algebraic analogue is
available for Frobenius Poisson algebras.  Recall that a finite-dimensional
Poisson algebra is Frobenius if its underlying algebra admits a
nondegenerate symmetric bilinear form invariant under multiplication;
equivalently, the induced map from the algebra to its linear dual is an
isomorphism of bimodules, up to the relevant grading shift.  Zhu, Van
Oystaeyen, and Zhang proved in \cite{ZVOZ} that the Poisson cochain
complexes with coefficients in a Frobenius Poisson algebra and in its
linear dual form a differential calculus.  In the present setting, the
differential calculus on polyvector fields and mixed differential forms
descends to Poisson cohomology and mixed Poisson homology.

\begin{theorem}[Theorems~\ref{DC-on-Poisson} and~\ref{mainthm22}]
\label{mainthm2}
(1)
Let $B=(\mathcal O(\Sigma),\pi)$ be a Poisson algebra with
$\deg(\pi)=(0,2)$.  Then the preceding construction endows
\[
\bigl(
\mathrm{HP}^{\bullet,\bullet}(B),
\widetilde{\mathrm{HP}}_{\bullet,\bullet}(B)
\bigr)
\]
with the structure of a differential calculus.

(2)
Let
$A=\bigl(\mathcal O(k^{m||n}),\pi\bigr)$,
$A^!=\bigl(\mathcal O(k^{n\wedge m}),\pi^!\bigr)$
be Koszul-dual admissible quadratic Poisson algebras.  Then the maps
induced by \eqref{eq:intro-calculus-correspondence} define an isomorphism
of differential calculi
\[
\phi_*:
\bigl(
\mathrm{HP}^{\bullet,\bullet}(A),
\widetilde{\mathrm{HP}}_{\bullet,\bullet}(A)
\bigr)
\xrightarrow{\;\sim\;}
\bigl(
\mathrm{HP}^{\bullet,\bullet}(A^!),
\widetilde{\mathrm{HP}}_{\bullet,\bullet}(A^!)
\bigr).
\]
\end{theorem}

Here $\mathrm{HP}^{\bullet,\bullet}(-)$ denotes Poisson cohomology; the
mixed Poisson homology
$\widetilde{\mathrm{HP}}_{\bullet,\bullet}(-)$ is defined in
Section~\ref{sect:Poisson}.

\subsection{Unimodularity
and the Batalin-Vilkovisky algebra}

For an oriented manifold, contraction with a volume form identifies its
space of polyvectors with its space of differential forms. If the
manifold is Poisson, this identification need not intertwine the
Poisson cochain and chain differentials, because a chosen volume form
need not be a Poisson cycle. The Poisson structure is called
\textit{unimodular} if there exists a volume form that is a Poisson
cycle; equivalently, its modular class vanishes. This notion was introduced by Weinstein (see \cite{Weinstein}).
For unimodular Poisson manifolds, Xu \cite{xu1999gerstenhaber}
established a duality between Poisson cohomology and Poisson
homology. For a related construction on Poisson manifolds with diagonalizable
modular vector fields, see \cite{chen2023batalin}. The analogous notion for Frobenius Poisson algebras is
discussed in \cite{ZVOZ}. 
 In both settings, a compatible volume form equips the
resulting structure with the structure of a {\it differential calculus with duality} in the sense of Lambre
\cite{lambre2010dualite}.

An important consequence of unimodularity, and more generally of a
differential calculus with duality, is that the corresponding cohomology
carries a Batalin--Vilkovisky algebra structure.
Such algebras encode the genus-zero algebraic structure of a topological conformal field theory; see,
for example, \cite{getzler1994batalin}.

It turns out that the above two differential
calculus structures are isomorphic
in the case of $k^{m|0}$
and $k^{0|m}$, when equipped
with Koszul dual quadratic Poisson
structures (see \cite{chen2021poisson}).
The following theorem extends this result to the mixed case.

\begin{theorem}[Theorem \ref{thm:UnimodtoBV} and \ref{mainthm33}]
\label{mainthm3}
(1)
Let $\Sigma$ be a finite-dimensional bigraded vector space, let $\pi$ be a Poisson bivector of bidegree $(0,2)$.  If $\pi$ is
unimodular with respect to a mixed volume form $\gamma$, then
$\bigl(
\mathrm{HP}^{\bullet,\bullet}(\mathcal O(\Sigma)),
\widetilde{\mathrm{HP}}_{\bullet,\bullet}(\mathcal O(\Sigma))
\bigr)$
forms a differential calculus with duality,
and therefore
$\mathrm{HP}^{\bullet,\bullet}(\mathcal O(\Sigma))$ is a Batalin--Vilkovisky algebra.

(2) Let
$A=(\mathcal O(k^{m||n}),\pi),
A^!=(\mathcal O(k^{n\wedge m}),\pi^!)$
be Koszul-dual admissible quadratic Poisson algebras, equipped with the mixed volume
forms in \eqref{volumeform} and \eqref{volumeform2}. Then $\pi$ is
unimodular if and only if $\pi^!$ is unimodular. In this case, the
induced map
\[
\phi_*:
\mathrm{HP}^{\bullet,\bullet}(A)
\xrightarrow{\;\sim\;}
\mathrm{HP}^{\bullet,\bullet}(A^!)
\]
is an isomorphism of Batalin--Vilkovisky algebras.
\end{theorem}

The paper is organized as follows. Section~\ref{sect:koszuldual} reviews Koszul duality for bigraded quadratic algebras and applies the general construction to the pair $\mathcal O(k^{m\Vert n})$ and $\mathcal O(k^{n\wedge m})$. Section~\ref{sect:diffcalcul} develops the differential calculus of polyvector fields and mixed differential forms on bigraded spaces and proves that the calculi associated with this Koszul-dual pair are isomorphic. Section~\ref{sect:Poisson} introduces the Poisson cochain and mixed Poisson chain complexes, proves that the resulting Poisson cohomology and mixed Poisson homology form a differential calculus, and shows that, in the unimodular case, a mixed volume form gives rise to a differential calculus with duality and hence to a Batalin--Vilkovisky algebra structure on Poisson cohomology. Section~\ref{sect:koszuldual5} studies admissible quadratic Poisson structures on the Koszul-dual pair and proves that Koszul duality identifies the associated differential calculi, preserves unimodularity, and induces an isomorphism of the resulting Batalin--Vilkovisky algebras.

\begin{remark}[Relation to previous work]
This paper overlaps in part with \cite{chen2021poisson},
where the Koszul duality of Poisson
structures on $k^{m|0}$
and $k^{0|m}$,
which corresponds to
$k^{m||0}$ and $k^{0\wedge m}$ in the current paper,
is studied.
In the mixed case, the commutation relations and differential calculus operations depend on both the parity degree and the weight. We formulate the isomorphisms of differential calculi using this bigrading.
Shoikhet observed in \cite{shoikhet2010Koszul} that, for
$k^{m|0}$
and $k^{0|m}$,
a bivector with quadratic
coefficients $\pi$
on $k^{m|0}$ is Poisson
if and only if the corresponding dual
bivector $\pi^!$ on $k^{0|m}$
is Poisson.
He then showed that, for $k=\mathbb{C}$, the Kontsevich deformation
quantizations of these two Poisson structures obtained through Tamarkin's approach \cite{Tamarkin} remain Koszul dual.
We extend the classical part of Shoikhet's results to bigraded spaces by
establishing the corresponding differential calculus isomorphisms and
the preservation of unimodularity and the induced Batalin--Vilkovisky
structures. Deformation quantization is pursued in subsequent work on two-brane formality with commutation factors.
\end{remark}

\section{Koszul duality for bigraded algebras}\label{sect:koszuldual}
This section reviews Koszul duality for {\it bigraded}
algebras. Koszul algebras were introduced by Priddy in \cite{Priddy}.
For a comprehensive treatment, see Loday--Vallette
\cite{loday2012algebraic}; our notation differs slightly from theirs.

\subsection{Bigraded quadratic and Koszul algebras}

Throughout this paper, all vector spaces are endowed with a 
$\mathbb{Z}\times\mathbb{Z}$-grading.
Recall that a
\textit{bigraded  vector space} $V$
is a direct sum of vector spaces indexed by $\mathbb{Z}\times\mathbb{Z}$; 
that is, $V=\bigoplus_{(i,j)\in \mathbb{Z}\times\mathbb{Z}}V_{i,j}$.
An element $a\in V_{i,j}$
is called {\it homogeneous of
degree} $(i,j)$. We call $i$ and $j$ the {\it parity} and 
{\it weight} of $a$,
and denote them by $p(a)$ and $w(a)$ 
respectively. The super parity of \(a\) is the reduction of \(p(a)\) modulo \(2\).
A linear map $f: V\to W$ between bigraded
vector spaces is said to have degree $(i,j)$ if for every $(s,t)\in \mathbb{Z}\times\mathbb{Z}$,
$f(V_{s,t})\subseteq W_{s+i,t+j}$.

\begin{convention}\label{conv:bigrading}
For a homogeneous element $a$, write $|a|=(p(a),w(a))$.  For homogeneous
elements $a$ and $b$, set
\[
|a||b|:=p(a)p(b)+w(a)w(b)\pmod 2.
\]
Interchanging $a$ and $b$ contributes the sign $(-1)^{|a||b|}$.
Tensor length, differential-form degree, and polyvector arity are recorded
separately and do not contribute to this sign.  We set $V[r]^{p,q}:=V^{p,q+r}$.  The suspension $s:V\to V[-1]$ has
bidegree $(0,1)$.
\end{convention}

\begin{remark}[Relation with graded super vector spaces]
The $\mathbb Z\times\mathbb Z$-grading used here refines the
$\mathbb Z_2\times\mathbb Z$-grading of the graded super vector
spaces considered by Deligne and Morgan
\cite{deligne1999notes}.
Indeed, a bigraded vector space
$V=\bigoplus_{i,j\in\mathbb Z}V_{i,j}$ determines a graded
super vector space by setting
\[
V^{\mathrm{super}}_{\epsilon,j}
:=
\bigoplus_{\substack{i\in\mathbb Z\\ i\bmod 2=\epsilon}}
V_{i,j},
\qquad
\epsilon\in\mathbb Z_2,\quad j\in\mathbb Z.
\]
The super parity of a homogeneous element $a$ is then
$p(a)\bmod 2$, while its weight remains $w(a)$.
The commutativity sign
$(-1)^{p(a)p(b)+w(a)w(b)}$ is unchanged by this reduction,
and the suspension still preserves super parity and raises
the weight by one.

We retain the integer parity degree to keep track of the full
bidegrees of homogeneous elements and maps.
\end{remark}

\begin{definition}
An associative algebra $A$ is called
a \textit{bigraded algebra} if it admits a 
direct sum $A=\bigoplus_{(i,j)\in \mathbb{Z}\times 
\mathbb{Z}}A_{i,j}$ of vector spaces with
a product that maps $A_{i_1,j_1}\times A_{i_2,j_2}$ to 
$A_{i_1+i_2,j_1+j_2}$.
\end{definition}

For a bigraded vector space $V$, its tensor algebra
$T(V)=\bigoplus_{r\geq0}V^{\otimes r}$ carries the internal bigrading and,
separately, the tensor length grading $r$.

Suppose that $V$ is a 
finite-dimensional
vector space. 
Let $R\subset V\otimes V$ be a 
subspace spanned 
by homogeneous elements.
We call the pair $(V, R)$ a \textit{quadratic datum}.

\begin{definition}
The {\it quadratic algebra} 
associated to the quadratic datum $(V, R)$ is $$A=A(V, R):=T(V)/(R),$$
where $(R)$ is
the two-sided ideal generated by $R$ in $T(V)$; that is,
\begin{equation*}
A=T(V)/(R)=k\oplus V\oplus 
\frac{V^{\otimes 2}}{R}\oplus  \cdots\oplus 
\frac{V^{\otimes n}}{\sum_{i=0}^{n-2} V^{\otimes i}
\otimes R\otimes V^{\otimes n-2-i}}\oplus \cdots.
\end{equation*}
\end{definition}

\begin{definition}
Suppose that $A$ is the quadratic
algebra associated to $(V, R)$.

(1) The \textit{quadratic dual algebra}
of $A$,
denoted by $A^!$,
is the quadratic algebra
associated to 
$(sV^*,s\otimes s
\circ R^{\perp})$, where 
$V^*:=\mathrm{Hom}(V,k)$ is the
$k$-linear dual of $V$,
$s$ is the degree shifting functor
which preserves the parity and shifts
the weight up by $1$,
and $R^\perp$ is the annihilator
of $R$. 

(2) The \textit{quadratic dual coalgebra} of \( A \),
denoted by \( A^\ac \), is
the tensor-length graded dual
\[
A^\ac:=(A^!)^{\#}:=\bigoplus_{r\geq0}(A^!_r)^*.
\]
Since $V$ is finite-dimensional, every tensor-length component is finite
dimensional, and this graded dual is well defined.
\end{definition}

It follows directly from the definition that
$$A^{\ac}=\bigoplus_{n=0}^\infty A_n^{\ac}
=\bigoplus_{n=0}^\infty\bigcap_{i=0}^{n-2}
(s^{-1}V)^{\otimes i}\otimes
(s^{-1}\otimes s^{-1}(R))
\otimes (s^{-1}V)^{n-i-2},
$$
where it is understood that $A_0^{\ac}=k$ and
$A_1^{\ac}=s^{-1}V$.
Thus $A^{\ac}$ inherits a
coalgebra structure from the deconcatenation coproduct on
$T(s^{-1}V)$.
Moreover,
one has $(A^!)^!=A$.

Let
\[
\kappa:A^\ac\twoheadrightarrow s^{-1}V\xrightarrow{s}V
\rightarrowtail A
\]
be the canonical twisting morphism.  If
$\Delta_{1,n-1}(c)=\sum c_{(1)}\otimes c_{(2)}$ denotes the component of
the deconcatenation coproduct whose first factor has tensor length one,
define
\[
d_\kappa(a\otimes c)
=\sum(-1)^{w(a)}a\kappa(c_{(1)})\otimes c_{(2)}.
\]
The sign is the braided tensor-product sign contributed by the bidegree
$(0,1)$ of $\kappa$.
The following sequence is a chain complex associated with $A$, called the 
\textit{left Koszul complex}:
\begin{align}\label{koszulcomplex}
\cdots \stackrel{d_\kappa}
{\longrightarrow}A 
\otimes A^\ac_{i+1} \stackrel{d_\kappa}{\longrightarrow} 
A \otimes A^\ac_i 
\stackrel{d_\kappa}{\longrightarrow} \cdots
\longrightarrow A \otimes A^\ac_0 
\stackrel{\epsilon}\longrightarrow k,
\end{align}
where $\epsilon:A\to k$ is the augmentation.  The quadratic relations
imply $d_\kappa^2=0$.  Interchanging the
tensor factors gives the analogous {\it right Koszul complex}.
The left Koszul complex of $A$ can be written as the left twisted tensor product 
$A\otimes_\kappa A^\ac$.
Similarly, the right Koszul complex of $A$ can be written as 
$A^\ac\otimes_\kappa A$.

 \begin{definition}\label{DefKoszul}
 Suppose that $A$ is a quadratic algebra.  Then $A$ is called a
 {\it Koszul algebra} if the augmented Koszul complex
 \eqref{koszulcomplex} is exact.
 \end{definition}

If $A=A(V,R)$ is Koszul, then $A^!$ is also Koszul. Indeed, \cite[Theorem 2.3.1]{loday2012algebraic} shows that acyclicity of the left 
Koszul complex is equivalent to acyclicity of the right Koszul complex. 
The Koszul complex of $A^!$ is $$\cdots 
\stackrel{d}\longrightarrow A^! \otimes 
(A^!)^\ac_{i+1}\stackrel{d}\longrightarrow A^! \otimes 
(A^!)^\ac_i\stackrel{d}\longrightarrow\cdots 
\longrightarrow A^! \otimes (A^!)^\ac_0 
\stackrel{d}\longrightarrow k,$$
and its tensor-length graded dual is the right Koszul complex of
$A$, up to the suspension fixed above.
If $A$ is Koszul, we henceforth call 
the quadratic dual algebra $A^!$ and coalgebra $A^{\ac}$
its \textit{Koszul dual algebra} and \textit{Koszul dual coalgebra},
respectively.

\subsection{The Koszul dual pair
$\mathcal O(k^{m||n})$ and $\mathcal O(k^{n\wedge m})$}

We focus on the following two bigraded algebras:
\begin{equation}\label{eq:newpolyRmn}
\mathcal{O} (k^{m||n}) := \mathrm{Sym}\{ x_1,\cdots,x_m, y_1,\cdots,y_n \},  
\end{equation}
with 
$\mathrm{deg}(x_i)=(0,0)$
and 
$\mathrm{deg}(y_j)= (1,0)$ and
\begin{equation}\label{eq:newpolyRnm}
\mathcal{O} (k^{n\wedge m} )
:= \mathrm{Sym}\{ \xi_1,\cdots,\xi_m, \eta_1,\cdots,\eta_n \},
\end{equation}
with 
$
\mathrm{deg}(\xi_i)=(0,1)$
and
$\mathrm{deg}(\eta_j)=(-1,1)$.
The first main result is the following.

\begin{theorem}[Theorem
\ref{mainthm1}]\label{newmainthm1}
Let $\mathcal{O}(k^{m||n})$
and $\mathcal{O}(k^{n\wedge m})$
be as above.
Then both are Koszul algebras, and they are Koszul dual to each other.
\end{theorem}

\begin{proof}
By \eqref{eq:newpolyRmn},
\(\mathcal{O}(k^{m||n})=T(V)/(R)\), where
\begin{equation}\label{quadratic}
    V=\mathrm{span}\{x_1,\dots,x_m,y_1,\dots,y_n\}
\end{equation}
and
\[
R=\mathrm{span}\left\{
x_ix_j-x_jx_i,\ y_uy_v+y_vy_u,\ x_iy_u-y_ux_i
\,\middle|\,
\substack{1\leq i,j\leq m\\1\leq u,v\leq n}
\right\}.
\]
By quadratic duality,
\[
\mathcal{O}(k^{m||n})^!
=T(sV^*)/(s\otimes s\circ R^\bot),\qquad
V^*=\mathrm{span}\{x_1^*,\dots,x_m^*,y_1^*,\dots,y_n^*\},
\]
where
\[
R^\bot=\mathrm{span}\left\{
x_i^*x_j^*+x_j^*x_i^*,\
y_u^*y_v^*-y_v^*y_u^*,\
x_i^*y_u^*+y_u^*x_i^*
\,\middle|\,
\substack{1\leq i,j\leq m\\1\leq u,v\leq n}
\right\}.
\]
On the other hand, \eqref{eq:newpolyRnm} gives
\(\mathcal{O}(k^{n\wedge m})=T(\tilde V)/(\tilde R)\), where
\(\tilde V=\mathrm{span}\{\xi_1,\dots,\xi_m,\eta_1,\dots,\eta_n\}\)
and
\[
\tilde R=\mathrm{span}\left\{
\xi_i\xi_j+\xi_j\xi_i,\
\eta_u\eta_v-\eta_v\eta_u,\
\xi_i\eta_u+\eta_u\xi_i
\,\middle|\,
\substack{1\leq i,j\leq m\\1\leq u,v\leq n}
\right\}.
\]
Under the identifications
\(sx_i^*\leftrightarrow\xi_i\) and
\(sy_j^*\leftrightarrow\eta_j\), for
\(i=1,\dots,m\) and \(j=1,\dots,n\), we have
\(\tilde V\simeq sV^*\) and
\(\tilde R\simeq s\otimes s\circ R^\bot\).
Thus \(\mathcal{O}(k^{n\wedge m})\cong
\mathcal{O}(k^{m||n})^!\), and the two algebras are quadratic dual.

We next prove that \(\mathcal{O}(k^{m||n})\) is Koszul.
For homogeneous elements
\(v_1\wedge\cdots\wedge v_q\in\mathcal{O}(k^{m||n})_q\) and
\(u_1\wedge\cdots\wedge u_r\in
{\mathcal{O}(k^{m||n})^\ac}_r\), interchanging two adjacent factors
produces the signs
\((-1)^{p(v_i)p(v_{i+1})+w(v_i)w(v_{i+1})}\) and
\((-1)^{p(u_i)p(u_{i+1})+w(u_i)w(u_{i+1})}\), respectively.
The Koszul boundary map and the homotopy are given by
\begin{align*}
d:\mathcal{O}(k^{m||n})_q\otimes
{\mathcal{O}(k^{m||n})^\ac}_r
&\longrightarrow
\mathcal{O}(k^{m||n})_{q+1}\otimes
{\mathcal{O}(k^{m||n})^\ac}_{r-1},\\
(v_1\wedge\cdots\wedge v_q,\,
u_1\wedge\cdots\wedge u_r)
&\longmapsto
\sum_{j=1}^{r}(-1)^\epsilon
\bigl(v_1\wedge\cdots\wedge v_q\wedge su_j,\,
u_1\wedge\cdots\wedge\widehat{u_j}\wedge\cdots\wedge u_r\bigr),\\
h:\mathcal{O}(k^{m||n})_q\otimes
{\mathcal{O}(k^{m||n})^\ac}_r
&\longrightarrow
\mathcal{O}(k^{m||n})_{q-1}\otimes
{\mathcal{O}(k^{m||n})^\ac}_{r+1},\\
(v_1\wedge\cdots\wedge v_q,\,
u_1\wedge\cdots\wedge u_r)
&\longmapsto
\sum_{j=1}^{q}(-1)^\varepsilon
\bigl(v_1\wedge\cdots\wedge\widehat{v_j}\wedge\cdots\wedge v_q,\,
 s^{-1}v_j\wedge u_1\wedge\cdots\wedge u_r\bigr),
\end{align*}
where
\(\epsilon=\sum_{s=1}^{j-1}
(p(u_j)p(u_s)+w(u_j)w(u_s))\) and
\(\varepsilon=\sum_{t=j+1}^{q}
(p(v_j)p(v_t)+w(v_j)w(v_t))\).
A direct computation gives
\(d\circ h+h\circ d=(r+q)\operatorname{id}\).  Since $k$ has
characteristic zero, $(r+q)^{-1}h$ is a contracting homotopy in every
positive total tensor length; the degree-zero term is the augmentation.
Hence the augmented Koszul complex is exact. Therefore,
\(\mathcal{O}(k^{m||n})\) is Koszul. Since
\(\mathcal{O}(k^{n\wedge m})\) is its quadratic dual,
\(\mathcal{O}(k^{n\wedge m})\) is also Koszul, and the two algebras
are Koszul dual.
\end{proof}

Recall that 
the algebra $\mathcal O(k^{m|n})$ 
of the graded space
$k^{m|n}$
(see, for example, \cite{deligne1999notes})
is a $\mathbb{Z}$-graded
commutative
algebra 
generated by
$\{\tilde{x}_1,\cdots,\tilde{x}_m,
\tilde{y}_1,\cdots,\tilde{y}_n\},$
where
$\deg(\tilde{x}_i)= 0$ and
$\deg(\tilde{y}_j) = 1$.
We regard this algebra as bigraded by setting
$\deg(\tilde x_i)=(0,0)$ and $\deg(\tilde y_j)=(1,0)$.
Consider the following map
of associative algebras
$$\mathcal{O}(k^{m|n})
\to
\mathcal{O}(k^{m||n}):
\tilde x_i\mapsto x_i,\,
\tilde y_j\mapsto y_j,
$$
for $i=1,\cdots, m,\  j=1,\cdots, n$.
It is straightforward to see that
this is an isomorphism
of bigraded commutative 
algebras. In summary, we have
the following.

\begin{proposition}\label{prop:Rm|n}
(1) There is a natural isomorphism of bigraded commutative algebras
\[
\mathcal O(k^{m|n})\cong\mathcal O(k^{m||n}).
\]

(2) In the two extreme cases $n=0$ and $m=0$, the underlying ordinary superspaces admit natural
identifications
\[
k^{m||0}\cong k^{m|0}\cong k^m,
\qquad
k^{0\wedge n}\cong k^{0|n}.
\]
\end{proposition}

Recall that the usual superspace $k^{m|m}$ is Koszul self-dual; however, in the bigraded case, $k^{m||m}$ is not isomorphic to $k^{m\wedge m}$ and hence are not Koszul self-dual.
\begin{remark}
In the rest of the paper, we consider the Koszul dual pair
$\mathcal O(k^{m||n})$ and $\mathcal O(k^{n\wedge m})$, where
$\mathcal O(k^{m||n})$ is obtained by equipping the coordinate algebra
of $k^{m|n}$ with the bidegrees specified above. Other choices of
bidegrees for the generators can be handled similarly.
\end{remark}

\section{Differential calculus structures in the bigraded setting}\label{sect:diffcalcul}
Building on the bigraded Koszul framework developed in Section~\ref{sect:koszuldual}, this section introduces the differential-geometric operations needed to formulate differential calculus in the bigraded setting. We first study differential forms, polyvector fields, and the Schouten--Nijenhuis bracket, and then construct contractions, mixed differential forms, and the corresponding Cartan identities. These constructions culminate in a differential calculus structure and prepare the passage to Poisson (co)homology in Section~\ref{sect:Poisson}.

\subsection{Differential forms and polyvector fields}

In this subsection, we consider differential forms and polyvector fields on a bigraded commutative algebra and introduce the Schouten--Nijenhuis bracket. This also fixes the notation and sign conventions used below.

\begin{definition}[K\"ahler differential]
Suppose $A$ is a bigraded commutative algebra.
The $A$-module of {\it K\"ahler differentials},
denoted by $\underline{\Omega}(A)$,
is the quotient of the free $A$-module
generated by symbols $\underline{d}x$ with
$\mathrm{deg}(\underline{d}x)=\mathrm{deg}(x)$, for all 
$x\in A$,
by the following relations:
$$\underline{d}(x+y)\sim 
\underline{d}x+\underline{d}y, \quad
\underline{d}
(xy)\sim(\underline{d}x) y+x
\underline{d}(y),\quad
\underline{d}a\sim 0,$$
where $x,y\in A$ and $a\in k$.
\end{definition}

The space of 1-forms $\Omega^1(A)$ on $A$ has the same underlying vector space as
$\underline{\Omega}(A)$ but with
the weights on the elements shifted by $-1$. 
The image of $\underline{d}x$, for any $x$,
under this weight shift is
denoted by $dx$.
Thus $d$ may be regarded as a $k$-linear operator
of degree $(0,-1)$
that vanishes on $k$ and
satisfies the bigraded Leibniz rule
$$d(xy)=(dx)y+(-1)^{w(x)}x dy.$$
In what follows, we shall use both
$\underline{\Omega}(A)$
and $\Omega^1(A)$.

More generally, we have the following.

\begin{definition}[Differential forms]
For $p\ge 0$, the $A$-module of {\it differential
$p$-forms}, denoted by $\Omega^p(A)$,
is the bigraded symmetric $A$-module generated by 
$\Omega^1(A)$; that is,
$$\Omega^p(A)=\mathrm{Sym}^p_A(\Omega^1(A)).$$
The {\it total
space of differential forms} of $A$ is
$$\Omega^\bullet(A)
:=\bigoplus_{p\ge 0}\Omega^p(A).$$
\end{definition}

Extend the degree $(0,-1)$ map
$$d: A\to\Omega^1(A),\quad a\mapsto da,$$
to $\Omega^\bullet(A)$ 
by derivation and by letting $d^2=0$.
Then $(\Omega^\bullet(A),d)$ 
is a {\it differential bigraded algebra}.

\begin{definition}[Derivation]
Suppose that $A$ is a bigraded commutative algebra and that $M$ is a
bigraded $A$-module.  A homogeneous $k$-linear map
$\gamma:A\to M$ of bidegree $|\gamma|$ is a \textit{derivation} if
\[
\gamma(ab)=\gamma(a)b+(-1)^{|\gamma||a|}a\gamma(b)
\]
for all homogeneous $a,b\in A$, and if $\gamma(c)=0$ for $c\in k$.
The bigraded vector space of derivations is denoted by
$\overline{\mathfrak{X}}(A,M)$.
If $M=A$, we abbreviate $\overline{\mathfrak{X}}(A,M)$ as $\overline{\mathfrak{X}}(A)$.
\end{definition}

Here and below, a homogeneous $A$-linear map $F:M\to N$ of bidegree
$|F|$ is:
\[
F(am)=(-1)^{|F||a|}aF(m)
\]
for homogeneous $a\in A$ and $m\in M$.

\begin{proposition}
For a bigraded commutative algebra $A$ and an
$A$-module $M$, there is an isomorphism
$$\overline{\mathfrak{X}}(A,M)\cong
\mathrm{Hom}_A(\underline{\Omega}(A),M).$$
\end{proposition}

\begin{proof}
For a homogeneous derivation $D$ of bidegree $|D|$, define
$\Phi(D)$ on generators by
\[
\Phi(D)(\underline d a)=D(a).
\]
The defining relation
$\underline d(ab)=(\underline d a)b+a\underline d b$ and braided
$A$-linearity give
\[
\begin{aligned}
\Phi(D)(\underline d(ab))
=D(a)b+(-1)^{|D||a|}aD(b)
=D(ab).
\end{aligned}
\]
Thus $\Phi(D)$ is well defined and has bidegree $|D|$.
Conversely, for a homogeneous braided $A$-linear map
$F:\underline\Omega(A)\to M$, set $D_F(a)=F(\underline d a)$.
The same calculation shows that $D_F$ is a derivation of bidegree
$|F|$.  The assignments $D\mapsto\Phi(D)$ and $F\mapsto D_F$ are
inverse to one another.
\end{proof}

Analogously, let
$\mathfrak{X}^1(A)$ be the $A$-module
which has the same underlying space as
$\overline{\mathfrak{X}}(A)$ but with the
weights of the elements shifted by $1$.

\begin{definition}[Polyvectors]
Suppose $A$ is a bigraded commutative algebra.
For $p\ge 0$, the $A$-module of
$p$-polyvectors, denoted by 
$\mathfrak{X}^p(A)$, is the $p$-th
bigraded symmetric power of 
$
\mathfrak{X}^1(A)
$; that is,
$$
\mathfrak{X}^p(A):=\mathrm{Sym}^p_A(\mathfrak{X}^1(A)).
$$
The space of polyvectors on $A$ is
$$\mathfrak{X}^\bullet(A):=
\bigoplus_{p\ge 0}\mathfrak{X}^p(A).$$
\end{definition}

In what follows, for two
vector fields $X$ and $Y$,
we usually write $X\wedge Y$
for their symmetric product.

\begin{example}
\label{ex:polyvec}
Let $A$ be the free bigraded commutative algebra generated by $\{z_1,\cdots,z_l\}$.
Then 
$$\mathfrak{X}^\bullet(A)=
\mathrm{Sym}^\bullet\{z_1,\cdots,z_l,
\partial_{z_1},\cdots,\partial_{z_{l}}\},$$
where $\partial_{z_i}$
corresponds,
under the degree shifting map
$\overline{\mathfrak{X}}(A)
\stackrel{\simeq}\to
\mathfrak{X}^1(A)$,
to the derivation 
$\overline{\partial}_{z_i}$
that maps $z_i$ to 1 and
all other $z_j$ to zero.

The degrees of the symbols appearing
in Example \ref{ex:polyvec}:

\begin{table}[h]
\centering
\begin{tabular}{ >{\centering\arraybackslash}
m{3cm} 
>{\centering\arraybackslash}
m{3cm} 
>{\centering\arraybackslash}
m{3cm} 
>{\centering\arraybackslash}
m{3.5cm} 
}
\hline
Variable 
 & Degree&Variable &Degree\\
\hline\hline
$\underline{d}x$
&$(p(x),w(x))$
&$dx$&
$(p(x),w(x)-1)$\\[2mm]
$\overline{\partial}_x$
&$(-p(x),-w(x))$
&$\partial_x$&
$(-p(x),-w(x)+1)$\\[2mm]
\hline
\end{tabular}
\caption{The degree assignments}
\label{table:0}
\end{table}
\end{example}

We next study the bigraded
Schouten-Nijenhuis bracket
defined on polyvectors.
For an element $D
\in\mathfrak{X}^1(A)$, we write
the corresponding derivation
in $\overline{\mathfrak{X}}(A)$
as $\overline{D}$, and conversely. We also use
$s:
\overline{\mathfrak{X}}(A)\to
\mathfrak{X}^1(A)$
for this identification.
To simplify notation, we write 
$$|x|:=\mathrm{deg}(x),|x|+|y|:=(p(x)+p(y),w(x)+w(y))
,|x||y|:=p(x)p(y)+w(x)w(y)\ \mathrm{and}\ \partial_i:=\partial_{z_i}.$$

Now suppose $A$ is a bigraded commutative
algebra. For
any $\overline{X},\overline{Y}
\in\overline{\mathfrak{X}}(A)$,
their commutator
$$[\overline{X}, \overline{Y}](f)
:=\overline{X}(\overline{Y}(f))-
(-1)^{|\overline{X}||\overline{Y}|}\overline{Y}
(\overline{X}(f)),\quad\mbox{for any}\,
f\in A,$$
defines a derivation of $A$.
\begin{definition}[Schouten-Nijenhuis
bracket]
For any two polyvectors
$Y_1\wedge\cdots\wedge Y_p$
and $Z_1\wedge\cdots\wedge Z_q$
in $\mathfrak{X}^\bullet(A)$,
their Schouten-Nijenhuis bracket
is given by
\begin{align}\label{SNbracket}
    &[Y_1\wedge\cdots \wedge Y_{p}\ , 
\  Z_{1}\wedge \cdots \wedge Z_{q}]\cr
&:=\sum^{p,q}_{s=1,t=1}  
(-1)^\varepsilon Y_{1}\wedge\cdots\wedge \widehat Y_{s}\wedge\cdots 
\wedge Y_{p}
\wedge s([\overline{Y}_{s},
\overline{Z}_{t}])\wedge Z_{1} \wedge \cdots \wedge \widehat Z_{t}\wedge\cdots
\wedge Z_{q},
\end{align}
where $\varepsilon=|Z_t|\sum_{i=1}^{t-1}|Z_i|+|Y_s|\sum_{j=s+1}^{p}|Y_j|$.
For homogeneous $f,g\in A$ and $X\in\mathfrak X^1(A)$,
set $[f,g]=0$ and $[X,f]=\overline X(f)$.
\end{definition}

The Schouten--Nijenhuis bracket
has the following properties:
\begin{enumerate}
\item[(1)] the bigraded skew commutativity:
\begin{align*}
[P,Q]=-(-1)^{(|P|-|s|)(|Q|-|s|)}[Q,P]; 
\end{align*}

\item[(2)] the bigraded Jacobi identity:
\begin{align*}
  (-1)^{(|P|-|s|)(|S|-|s|)}\big[[P,Q],S\big] 
  + (-1)^{(|P|-|s|)(|Q|-|s|)}\big[[Q,S],P\big]\\ 
   + (-1)^{(|Q|-|s|)(|S|-|s|)}\big[[S,P],Q\big] = 0; 
\end{align*}

\item[(3)] the bigraded Leibniz rule:  
\begin{align*}
  [P,\, Q\wedge S] = [P,Q]\wedge S 
  + (-1)^{|Q|(|P|-|s|)} Q\wedge[P,S], 
\end{align*}
for any homogeneous elements $P,Q,S\in \mathfrak{X}^\bullet(A)$.
\end{enumerate}

\begin{example}[Example 
\ref{ex:polyvec} continued]\label{ex:SNB}
Let $\Sigma$ be a bigraded vector space 
with coordinates $\{z_1,\cdots,z_l\}$, and set $A=\mathcal{O}(\Sigma)$. 
Suppose the homogeneous vector fields $X$ and $Y$ are given by 
$X=f_1\partial_{z_1}+\cdots+f_l\partial_{z_l}$ and 
$Y=g_1\partial_{z_1}+\cdots+g_l\partial_{z_l}$.
Then their Schouten-Nijenhuis bracket is  
$$[X,Y]=\sum_{i=1}^{l}\sum_{j=1}^{l}
(f_j\overline{\partial}_{z_j}g_i-(-1)^{(|\overline{\partial}_j|+|g_j|)
(|\overline{\partial}_i|+|f_i|)}g_j\overline{\partial}_{z_j}f_i)\partial_{i}.
$$
Suppose
$P\in \mathfrak{X}^p(A)$ and $Q\in\mathfrak{X}^q(A)$
are two polyvectors on $\Sigma$, which are given by 
\begin{align*}
P=\sum_{i_1,\cdots,i_p}P^{i_1,\cdots,i_p}\partial_{{i_1}}
\wedge\cdots\wedge\partial_{{i_p}}, \quad
Q=\sum_{j_1,\cdots,j_q}Q^{j_1,\cdots,j_q}
\partial_{{j_1}}\wedge\cdots\wedge\partial_{{j_q}},
\end{align*}
where $P^{i_1,\cdots,i_p},Q^{j_1,\cdots,j_q}\in A$.
Their bracket is obtained from \eqref{SNbracket} by viewing each
coefficient times a coordinate vector field as a homogeneous vector field
and then applying the bigraded Leibniz rule.  This prescription fixes all
coordinate signs through the sign rule in Convention~\ref{conv:bigrading}
and does not depend on a chosen ordering of the displayed monomials.
\end{example}

\subsection{Contractions, mixed differential forms and mixed Cartan identities}
In this and the following subsections, we develop the differential calculus on the spaces of polyvectors and mixed differential forms. We begin by defining the contraction action and proving the bigraded mixed Cartan identity.

The interior product is the contraction 
of a differential form with a vector field.
\begin{definition}[Contraction]\label{Contraction}
For any vector field $X\in \mathfrak{X}^1(A)$, the contraction  
\begin{align*}
        \iota_X:\Omega^p(A)\to\Omega^{p-1}(A), p\ge 1,
\end{align*}
is the $A$-linear map such that
\begin{align*}
        (\iota_X\omega)(X_1,\cdots,X_{p-1})=\omega(X,X_1,\cdots,X_{p-1})
\end{align*}
for any vector fields $X_1,\cdots,X_{p-1}$.
When $\omega$ is a $0$-form, i.e., $\omega\in A$, $\iota_X\omega$ is set to $0$.
\end{definition}
The contraction map satisfies the Leibniz rule: 
for any homogeneous elements $\alpha,\beta\in \Omega^\bullet(A)$,
\begin{align}\label{contraction:Leibniz}
    \iota_X(\alpha\wedge\beta)=(\iota_X\alpha)\wedge\beta+
    (-1)^{|\alpha||X|}\alpha\wedge(\iota_X\beta).
\end{align}

\begin{example}[Example
\ref{ex:SNB} continued]
    In the coordinate system $\{z_1,\cdots,z_l\}$ of $\Sigma$, 
    a vector field has the form $X=f_1\partial_1+\cdots+
    f_l\partial_l\in \mathfrak{X}^1(A)$. The contraction map $\iota_X$ is given by 
\begin{align*}
    \iota_X(dz_{i_1}\wedge\cdots\wedge d{z_{i_p}})=
    \sum_{j=1}^{p}(-1)^{|dz_{i_j}|(|dz_{i_1}|+\cdots+|dz_{i_{j-
    1}}|)}f_{i_j}dz_{i_1}\wedge\cdots\wedge \widehat{dz_{i_j}}
    \wedge\cdots\wedge d{z_{i_p}}.
\end{align*}
\end{example}

\begin{definition}[Lie derivative]\label{Lie derviative}
\begin{enumerate}
\item[(1)]  For any two vector fields $X,Y$, the \textit{Lie derivative} of 
$Y$ with respect to $X$, 
denoted by $\mathcal{L}_XY$, is the Schouten-Nijenhuis bracket of $X$ and $Y$; that is, $\mathcal{L}_XY=[X,Y]$.

\item[(2)] The 
\textit{Lie derivative} of 
a differential form $\omega$
with
respect to $X$ is given
by $$\mathcal{L}_X\omega:=
[\iota_X,d](\omega)
=\iota_Xd\omega
-(-1)^{|d||X|} d(\iota_X
\omega).$$
\end{enumerate}
\end{definition}

\begin{proposition}\label{prop:cartan-identities}
The following properties hold:
\begin{enumerate}
\item [$(1)$] The Lie derivative of a function $f$ with respect to 
$X$ is $\overline{X}(f)$.
        
\item [$(2)$] The Lie derivative obeys the Leibniz rule: 
for any homogeneous differential forms $\alpha,\beta$,
$$\mathcal{L}_X(\alpha\wedge\beta)=\mathcal{L}_X(\alpha)\wedge\beta+
(-1)^{|\alpha||\mathcal{L}_X|}\alpha\wedge\mathcal{L}_X(\beta).$$

\item[$(3)$] Cartan's formula holds:
$[\mathcal{L}_X,\iota_Y]=\iota_{[X,Y]}$, for any homogeneous vector fields $X$ and $Y$.

\item[$(4)$] The Lie derivative commutes with the differential $d$ on differential forms: 
$[\mathcal{L}_X,d]=0.$ 
\end{enumerate}
\end{proposition}

\begin{proof}
(1) For any function $f$, $$\mathcal{L}_Xf=[\iota_X,d](f)=\iota_X(df)=\langle 
df,X\rangle=\overline{X}(f).$$

(2) For any homogeneous differential forms $\alpha,\beta$, 
we have 
\begin{align*}
[\iota_X,d](\alpha\wedge\beta)=&(\iota_X\circ d-(-1)^{|d||X|}d\circ\iota_X)(\alpha\wedge\beta)\\
=&\iota_X(d\alpha\wedge\beta+(-1)^{|d||\alpha|}\alpha\wedge d\beta)-
(-1)^{|d||X|}d(\iota_X\alpha\wedge\beta+(-1)^{|X||\alpha|}\alpha\wedge\iota_X\beta)\\
=&\iota_X(d\alpha)\wedge\beta+(-1)^{|d||\alpha|+|X||\alpha|}\alpha\wedge\iota_X(d\beta)\\
&-(-1)^{|d||X|}(d(\iota_X\alpha)\wedge\beta+
(-1)^{|X||\alpha|+|d||\alpha|}\alpha\wedge d(\iota_X\beta))\\
=&([\iota_X,d]\alpha)\wedge\beta+
(-1)^{|\alpha||X|+|d||\alpha|}\alpha\wedge([\iota_X,d]\beta).
\end{align*}

(3) For any $\omega\in \Omega^p(A)$, $p\ge0$, and $X,Y\in \mathfrak{X}^1(A)$,
we prove that $[[\iota_X,d],\iota_Y](\omega)=\iota_{[X,Y]}\omega$.
We proceed by induction on $p$.
When $p=0$ and $\omega\in \Omega^0(A)$, $[[\iota_X,d],\iota_Y]
(\omega)=0=\iota_{[X,Y]}\omega.$
Then we verify $[[\iota_X,d],\iota_Y](df)=\iota_{[X,Y]}(df)$ for any $f\in A$.
By the definition of the Schouten--Nijenhuis bracket, $\iota_{[X,Y]}(df)=[\bar{X},\bar{Y}](f)=\bar{X}(\bar{Y}(f))-(-1)^{|\bar{X}||\bar{Y}|}\bar{Y}(\bar{X}(f)).$
And \begin{align}\label{cartan formula df}
    [[\iota_X,d],\iota_Y](df)=&(\iota_X\circ d \circ \iota_Y-(-1)^{|d||X|}d\circ\iota_X\circ\iota_Y
-(-1)^{(|X|+|d|)|Y|}\iota_Y\circ \iota_X\circ d
\\
&+(-1)^{(|X|+|d|)|Y|+|d||X|}\iota_Y\circ d \circ \iota_X)(df)\nonumber\\
=&(\iota_X\circ d \circ \iota_Y
+(-1)^{(|X|+|d|)|Y|+|d||X|}\iota_Y\circ d \circ \iota_X)(df)\nonumber\\
=&\bar{X}(\bar{Y}(f))-(-1)^{|\bar X||\bar Y|}\bar{Y}(\bar{X}(f))\nonumber\\
=&\iota_{[X,Y]}(df).\nonumber
\end{align}

Now assume that the identity holds for $\omega\in \Omega^p(A)$, $p>0$. Consider $\iota_{[X,Y]}$ acting on $df\wedge\omega$, $f\in \Omega^0(A)$.
By \eqref{contraction:Leibniz}, \eqref{cartan formula df}, and the 
induction hypothesis, we have
\begin{align*}
\iota_{[X,Y]}(df\wedge\omega)&= \iota_{[X,Y]}
(df)\wedge\omega +(-1)^{|[X,Y]||df|} 
df\wedge\iota_{[X,Y]}(\omega)\\
&=[[\iota_X,d],\iota_Y](df)\wedge\omega+
(-1)^{|[X,Y]||df|}df\wedge[[\iota_X,d],\iota_Y](\omega)\\
&=[[\iota_X,d],\iota_Y](df\wedge\omega).
\end{align*}
This completes the induction.

(4) By direct computation, we have
\begin{align*}
[\mathcal{L}_X,d]&=[[\iota_X, d],d]\\ &= [\iota_X, d] \circ d -(-1)^{|d|(|X|+|d|)} d \circ [\iota_X, d] \\  
&=\iota_X \circ d\circ d -(-1)^{|X||d|} d \circ \iota_X\circ d+(-1)^{|X||d|}d\circ 
\iota_X \circ d -d\circ d \circ \iota_X\\ 
    &= 0.
\end{align*}
\end{proof}

We now introduce mixed differential forms on a bigraded space $\Sigma$ such that $A=\mathcal{O}(\Sigma)=\mathrm{Sym}\{z_1,\cdots,z_l\}$ with generators $\{z_1,\cdots,z_l\}$, where
$p(z_i)+w(z_i)$ is even for $i\le m$ and $p(z_j)+w(z_j)$ is odd for $j\ge m+1$.

Define from $A$ the bigraded algebra
\begin{align*}
   \tilde{\Omega}^0(A):= \mathrm{Sym}\{z_1,\cdots,z_m,
   \partial_{m+1},\cdots,\partial_l\}, 
\end{align*} which is a subalgebra of $\mathfrak{X}^\bullet(A)$.
Let $\tilde{d}$ be a derivation of degree $(0,-1)$ on $\tilde{\Omega}^0(A)$, 
which acts on generators by $$\tilde{d}(z_i) = dz_i, i\le m\quad \mathrm{and}
\quad  \tilde{d}(\partial_j) = z_j^*, j\ge m+1.$$
Then \begin{align*}
    \tilde{d}:\tilde{\Omega}^0(A)\to \tilde{\Omega}^1(A)
    :=\Big\{\sum_{i=1}^{m}f_idz_i+\sum_{j=m+1}^{l}f_jz_j^*\mid f_i,f_j\in \tilde{\Omega}^0(A)\Big\}.
\end{align*}
\begin{definition}[Mixed differential form]
    For $p\ge 0$, the $\tilde{\Omega}^0(A)$-module of 
    \textit{mixed differential $p$-forms}, denoted by $\tilde{\Omega}^p(A)$, 
    is a bigraded symmetric $\tilde{\Omega}^0(A)$-module generated by 
    $\tilde{\Omega}^1(A)$; that is, 
    \begin{align*}
        \tilde{\Omega}^p(A)=\mathrm{Sym}^p_{\tilde{\Omega}^0(A)}(\tilde{\Omega}^1(A)).
    \end{align*}
    The \textit{total space of mixed differential forms} of $A$ is
    \begin{align*}
        \tilde{\Omega}^\bullet(A):=\bigoplus_{p\ge 0}\tilde{\Omega}^p(A)
        =\mathrm{Sym}\{z_1,\cdots,z_m,\partial_{m+1},\cdots,\partial_l,dz_1,
        \cdots,dz_m,z_{m+1}^*,\cdots,z_l^*\}.
    \end{align*}
\end{definition}
Extend the degree $(0,-1)$ map $\tilde{d}$ to $\tilde{\Omega}^\bullet(A)$ 
by derivation and by letting $\tilde{d}^2=0$.
The map $\tilde{d}$ is a differential of degree $(0,-1)$, then we have
\begin{align*}
    \tilde{d}(\omega\wedge \eta)=\tilde{d}\omega\wedge\eta+
    (-1)^{|\omega||\tilde{d}|}\omega\wedge \tilde{d}\eta,
\end{align*}
for any $\omega, \eta\in \tilde{\Omega}^{\bullet}(A)$.

Let $\langle-,-\rangle$ be the pairing between the generators of 
$\mathfrak{X}^\bullet(A)$ and $\tilde{\Omega}^\bullet(A)$: 
$$\langle dz_i,\partial_{j}\rangle
=\delta_{ij}, i\le m\quad \mathrm{and}\quad \langle z_s^*,z_t\rangle
=(-1)^{w(z_s)}\delta_{st}, s\ge m+1.$$

\begin{definition}[Mixed vector field]
Denote the module of \textit{mixed vector fields} on $A$ by $\tilde{\mathfrak{X}}^1(A)$,  $$\tilde{\mathfrak{X}}^1(A):
=\Big\{\sum_{i=1}^{m}g_i\partial_{i}+\sum_{j=m+1}^{l}g_jz_j\mid g_i,g_j\in 
\tilde{\Omega}^0(A)\Big\}\subset {\mathfrak{X}}^\bullet(A).$$
For any $p\ge 1$, the space of \textit{mixed $p$-vector fields}, denoted by $\tilde{\mathfrak{X}}^p(A)$, is the bigraded symmetric $\tilde{\Omega}^0(A)$-module generated by $\tilde{\mathfrak{X}}^1(A)$; that is, $$\tilde{\mathfrak{X}}^p(A)
:=\mathrm{Sym}_{\tilde{\Omega}^0(A)}^p(\tilde{\mathfrak{X}}^1(A))\\\subset 
\mathfrak{X}^\bullet(A).$$
\end{definition}

A vector $X\in \tilde{\mathfrak{X}}^1(A)$ defines a map 
$\tilde{\iota}_X$ on $\tilde{\Omega}^\bullet(A)$ of the same degree as $X$.
On a mixed 1-form, the operator $\tilde{\iota}_X$ is given by
$\tilde{\iota}_X\omega=\omega(X)=\langle\omega,X\rangle$.
Then it induces a natural action of $\tilde{\mathfrak{X}}^1(A)$ on 
$\tilde{\Omega}^\bullet(A)$.

\begin{definition}[Mixed contraction]
For any vector $X\in \tilde{\mathfrak{X}}^1(A)$,  
    \begin{align*}
        \tilde{\iota}_X:\tilde{\Omega}^p(A)\to\tilde{\Omega}^{p-1}(A)
    \end{align*}
    is the $\tilde{\Omega}^0(A)$-linear map defined by the property that
    \begin{align}\label{mixedcontraction}
        (\tilde{\iota}_X\omega)(X_1,\cdots,X_{p-1})=\omega(X,X_1,\cdots,X_{p-1})
    \end{align}
    for any vectors $X_1,\cdots,X_{p-1}\in \tilde{\mathfrak{X}}^1(A)$.
    When $\omega$ is a mixed $0$-form, $\iota_X\omega=0$ by convention.
    The differential form $\tilde{\iota}_X\omega$ is called the 
    \textit{mixed contraction} of $\omega$ with $X$.
\end{definition}
The mixed contraction map satisfies the Leibniz rule: for any homogeneous 
elements $\alpha,\beta\in \tilde{\Omega}^\bullet(A)$, we have 
\begin{align}\label{mixedcontraction:Leibniz}
    \tilde{\iota}_X(\alpha\wedge\beta)=(\tilde{\iota}_X\alpha)\wedge\beta+
    (-1)^{|\alpha||X|}\alpha\wedge(\tilde{\iota}_X\beta).
\end{align}
For $b\in\widetilde\Omega^0(A)$, set $\tilde\iota_b(\omega)=b\omega$.
We extend the mixed contraction to $\mathfrak{X}^\bullet(A)$ by 
    $\tilde{\iota}_{X\wedge Y}:=\tilde{\iota}_X\circ\tilde{\iota}_Y.$

\begin{proposition}[Mixed Cartan's formula]\label{mixedcartan}
    Given a bigraded space $\Sigma$ and $A=\mathcal{O}(\Sigma)$ as above, 
    for any homogeneous polyvector field $P,Q$, 
    \begin{align}\label{mixedcartanformula}
[[\tilde{\iota}_P,\tilde{d}],\tilde{\iota}_Q]
=\tilde{\iota}_{[P,Q]}.
\end{align}
    \end{proposition}
\begin{proof}
The identification with shifted derivations of $\widetilde\Omega^0(A)$ is
$$
\partial_i\longmapsto s(\overline\partial_{z_i}),\quad i\le m,
\qquad
z_j\longmapsto(-1)^{w(z_j)}s(\overline\partial_{\partial_j}),
\quad j\ge m+1.
$$

Under the preceding identifications, the required identity is the bigraded
Cartan identity on the free bigraded commutative algebra
$\widetilde\Omega^0(A)$.  For mixed vector fields it is
Proposition~\ref{prop:cartan-identities}, applied to
$\widetilde\Omega^0(A)$.
For mixed functions, the identity follows from the graded Leibniz rule.
For decomposable polyvectors, the identity follows inductively from the
Leibniz rule for the Schouten--Nijenhuis bracket and from
$\widetilde\iota_{P\wedge Q}=\widetilde\iota_P\widetilde\iota_Q$.
Both sides are linear in $P$ and $Q$, so the identity holds for arbitrary
mixed polyvectors.
\end{proof}

\subsection{Differential calculus and Koszul duality}
With the mixed Cartan identity established, we now recast the preceding constructions in the framework of differential calculus introduced by Tamarkin and Tsygan \cite{tamarkin2005ring}. We prove that, for the algebra of functions on a bigraded vector space, polyvector fields and mixed differential forms, equipped with the operations defined above, form a differential calculus. Specializing to $A=\mathcal{O}(k^{m||n})$ and $A^!=\mathcal{O}(k^{n\wedge m})$, we further show that the correspondence $\phi$ intertwines the calculus operations, yielding the basic identification used in the Poisson setting of Section~\ref{sect:Poisson}.

\begin{definition}[Tamarkin-Tsygan]\label{diffcalculus}
Let $\mathrm H^{\bullet,\bullet}$ and $\mathrm M^{\bullet,\bullet}$
be bigraded vector spaces, and set $\delta=(0,-1)$.  A
\textit{differential calculus} is a sextuple
$$
(\mathrm{H}^{\bullet,\bullet},
\mathrm M^{\bullet,\bullet}, \cup, \iota, [-,-], d)
$$
satisfying the following conditions:
\begin{enumerate}
\item[(1)] $(\mathrm H^{\bullet,\bullet},\cup)$ is a bigraded
commutative algebra, and the bracket has bidegree $\delta$.  Equivalently,
the shifted space $\mathrm H[1]$ is a bigraded Lie algebra.  The
bracket is a biderivation:
\[
[P,Q\cup R]=[P,Q]\cup R
+(-1)^{(|P|+\delta)|Q|}Q\cup[P,R].
\]

\item[(2)] $\mathrm M^{\bullet,\bullet}$ is a module over
$(\mathrm H^{\bullet,\bullet},\cup)$ through contractions.  For
homogeneous $P$ of bidegree $|P|$, the operator $\iota_P$ has bidegree
$|P|$, and $\iota_{P\cup Q}=\iota_P\iota_Q$.

\item[(3)] The operator $d$ on $\mathrm M^{\bullet,\bullet}$ has
bidegree $\delta$ and satisfies $d^2=0$.  If
$L_P=[\iota_P,d]$, then
\begin{equation}\label{calculusformula}
\iota_{[P,Q]}=[L_P,\iota_Q]
=L_P\iota_Q-(-1)^{(|P|+\delta)|Q|}\iota_QL_P.
\end{equation}
\end{enumerate}
\end{definition}

In the following, if $\cup$, $\iota$, $[-,-]$ and $d$
are clear from context, we abbreviate a differential calculus as
$(\mathrm H^{\bullet,\bullet},
\mathrm M^{\bullet,\bullet})$.

For polyvector fields and mixed differential forms, we introduce a further direct-sum decomposition determined by the bidegree:
\[
\mathfrak{X}^{i,j}(A)
:=
\bigl\{
P\in \mathfrak{X}^{\bullet}(A)
\mid \deg(P)=(i,j)
\bigr\},
\qquad
\widetilde{\Omega}^{i,j}(A)
:=
\bigl\{
P\in \widetilde{\Omega}^{\bullet}(A)
\mid \deg(P)=(i,j)
\bigr\}.
\]
Under the above bigrading, we proceed to study the differential calculus structure introduced earlier.

\begin{theorem}\label{DC-on-bigraded}
Let $A=\mathcal{O}(\Sigma)$ be an algebra of algebraic functions of a bigraded vector space $\Sigma$. Then
\[
\bigl(
\mathfrak{X}^{\bullet,\bullet}(A),
\tilde\Omega^{\bullet,\bullet}(A),
\wedge,\tilde\iota,[-,-],\tilde d
\ \bigr)
\]
forms a differential calculus.
\end{theorem}

\begin{proof}
The wedge product is bigraded commutative.  By its construction from the
Schouten--Nijenhuis bracket on
$\mathfrak X^\bullet(\widetilde\Omega^0(A))$, the bracket has bidegree
$\delta=(0,-1)$ and satisfies shifted skew-symmetry, the shifted Jacobi
identity, and the biderivation identity.  Hence
$\mathfrak X^{\bullet,\bullet}(A)$ is a bigraded Gerstenhaber algebra.

The contraction defined in \eqref{mixedcontraction} has bidegree $|P|$:
\[
\widetilde\iota:
\mathfrak X^{s,t}(A)\otimes\widetilde\Omega^{u,v}(A)
\longrightarrow\widetilde\Omega^{u+s,v+t}(A).
\]
Its extension to decomposable polyvectors satisfies
$\widetilde\iota_{P\wedge Q}=\widetilde\iota_P\widetilde\iota_Q$.
Finally, $\widetilde d$ has bidegree $\delta$ and squares to zero.  The
remaining calculus identity is precisely Proposition~\ref{mixedcartan}:
\[
\tilde\iota_{[P,Q]}=[L_P,\tilde\iota_Q] :={L}_{P}\circ
\tilde\iota_Q-(-1)^{|L_P||Q|}\tilde\iota_Q\circ {L}_P.
\qedhere
\]
\end{proof}

We now specialize Theorem \ref{DC-on-bigraded} to the algebras \(A=\mathcal{O}(k^{m||n})\) and \(A^!=\mathcal{O}(k^{n\wedge m})\) considered here. Since both \(A\) and \(A^!\) satisfy its assumptions, the corresponding structures
\[
\bigl(\mathfrak{X}^{\bullet,\bullet}(A),\widetilde{\Omega}^{\bullet,\bullet}(A),
\wedge,\widetilde{\iota},[-,-],\widetilde{d}\ \bigr)
\quad\text{and}\quad
\bigl(\mathfrak{X}^{\bullet,\bullet}(A^!),\widetilde{\Omega}^{\bullet,\bullet}(A^!),
\wedge,\widetilde{\iota},[-,-],\widetilde{d}\ \bigr)
\]
form differential calculi.

We introduce the following correspondence \(\phi\), which will be used repeatedly throughout the remainder of the paper. 

\begin{table}[h]
\centering
\begin{tabular}{ >{\centering\arraybackslash}
m{4cm} 
>{\centering\arraybackslash}
m{4cm} 
>{\centering\arraybackslash}
m{3cm}}
\hline
Variables of $k^{m||n}$ &
Variables of $k^{n\wedge m}$ & Degree\\
\hline\hline
$x_i$ &$\displaystyle\partial_{\xi_i}$ & $(0,0)$\\[2mm]
$\displaystyle\partial_{x_i}$ &$-\xi_i$&$(0,1)$\\[2mm]
$dx_i$&$\xi_i^*$ &$(0,-1)$\\[2mm]
$y_j$ & $\displaystyle\partial_{\eta_j}$ &$(1,0)$\\[2mm]
$\displaystyle\partial_{y_j}$& $\eta_j$ &$(-1,1)$\\[2mm]
$y_j^*$& $d\eta_j$ &$(-1,0)$\\[2mm]
\hline
\end{tabular}
\caption{The correspondence \(\phi\) between
$\mathcal{O}(k^{m||n})$
and $\mathcal{O}(k^{n\wedge m})$}
\label{table:1}
\end{table}
The correspondence $\phi$ in Table~\ref{table:1} induces two
bidegree-preserving maps:
\[
\phi_{\mathfrak X}:
\mathfrak X^{\bullet,\bullet}(A)
\longrightarrow
\mathfrak X^{\bullet,\bullet}(A^!),
\qquad
\phi_{\Omega}:
\widetilde\Omega^{\bullet,\bullet}(A)
\longrightarrow
\widetilde\Omega^{\bullet,\bullet}(A^!).
\]
For simplicity, we write
$\phi=(\phi_{\mathfrak X},\phi_{\Omega})$.

\begin{proposition}\label{prop:fourier-compatibility}
For homogeneous
$P,Q\in\mathfrak X^{\bullet,\bullet}(A)$ and
$\omega\in\widetilde\Omega^{\bullet,\bullet}(A)$,
\begin{align}
\phi_{\mathfrak X}(P\wedge Q)
 &=\phi_{\mathfrak X}(P)\wedge\phi_{\mathfrak X}(Q),
 \label{eq:phi-product}\\
\phi_{\mathfrak X}([P,Q])
 &=[\phi_{\mathfrak X}(P),\phi_{\mathfrak X}(Q)],
 \label{eq:phi-bracket}\\
\phi_\Omega(\widetilde d\omega)
 &=\widetilde d\,\phi_\Omega(\omega),
 \label{eq:phi-differential}\\
\phi_\Omega(\widetilde\iota_P\omega)
 &=\widetilde\iota_{\phi_{\mathfrak X}(P)}\phi_\Omega(\omega).
 \label{eq:phi-contraction}
\end{align}
\end{proposition}

\begin{proof}
In the present case, the polarized coordinate algebras are
\[
\widetilde\Omega^0(A)=\mathrm{Sym}\{x_1,\ldots,x_m,
\partial_{y_1},\ldots,\partial_{y_n}\},
\qquad
\widetilde\Omega^0(A^!)=\mathrm{Sym}\{\partial_{\xi_1},\ldots,
\partial_{\xi_m},\eta_1,\ldots,\eta_n\}.
\]
The degree-preserving assignments
$x_i\mapsto\partial_{\xi_i}$ and
$\partial_{y_j}\mapsto\eta_j$ extend uniquely to an isomorphism
$F:\widetilde\Omega^0(A)\to\widetilde\Omega^0(A^!)$ of free bigraded
commutative algebras.  Functoriality
of K\"ahler differentials and derivations sends
\[
dx_i\mapsto\xi_i^*,\qquad y_j^*\mapsto d\eta_j,
\qquad
\partial_{x_i}\mapsto-\xi_i,\qquad y_j\mapsto\partial_{\eta_j}.
\]
These are exactly the remaining rows of Table~\ref{table:1}.
The de Rham differential, contraction, wedge product, and
Schouten--Nijenhuis bracket are functorial under $F$.  Transporting them
through the identifications with the two mixed calculi gives
\eqref{eq:phi-product}--\eqref{eq:phi-contraction}.
\end{proof}

\begin{theorem}\label{thm:dc-koszul-duality}
Let \(A=\mathcal O(k^{m||n})\) and
\(A^!=\mathcal O(k^{n\wedge m})\).
The correspondence \(\phi\) in Table~\ref{table:1} induces an isomorphism
\[
\bigl(
\mathfrak X^{\bullet,\bullet}(A),
\widetilde\Omega^{\bullet,\bullet}(A),
\wedge,\widetilde\iota,[-,-],\widetilde d
\ \bigr)
\cong
\bigl(
\mathfrak X^{\bullet,\bullet}(A^!),
\widetilde\Omega^{\bullet,\bullet}(A^!),
\wedge,\widetilde\iota,[-,-],\widetilde d
\ \bigr)
\]
of differential calculi.
\end{theorem}

\begin{proof}
From the above definition, we have
\begin{eqnarray*}
    \mathfrak{X}^{\bullet,\bullet}(A)&=&\mathrm{Sym}^{\bullet,\bullet}\{x_1,\cdots,x_m,\partial_{y_1},\cdots,\partial_{y_n},\partial_{ x_1},\cdots,\partial_{x_m},y_1,\cdots,y_n\}
\end{eqnarray*} and \begin{eqnarray*}
    \mathfrak{X}^{\bullet,\bullet}(A^!)&=&\mathrm{Sym}^{\bullet,\bullet}\{\partial _{\xi_1},\cdots,\partial _{\xi_m},\eta_1,\cdots,\eta_n,\xi_1,\cdots,\xi_m,\partial_{ \eta_1},\cdots,\partial_{\eta_n}\},
\end{eqnarray*}
where $\mathrm{Sym}^{i,j}\{z_1,\cdots,z_l\}:=\{f\in\mathrm{Sym}\{z_1,\cdots,z_l\}\mid\mathrm{deg}(f)=(i,j) \}$.
Under the correspondence given in Table \ref{table:1},
we obtain an isomorphism of bigraded spaces 
$
\phi:\mathfrak{X}^{\bullet,\bullet}(A)\cong\mathfrak{X}^{\bullet,\bullet}(A^!).
$

The two maps are bijective because Table~\ref{table:1} gives mutually
inverse correspondences between homogeneous algebra generators.
Equations~\eqref{eq:phi-product} and~\eqref{eq:phi-bracket} identify the
Gerstenhaber products and brackets, while
\eqref{eq:phi-differential} and~\eqref{eq:phi-contraction} identify the
module differential and contraction action.  These are exactly the
operations in Definition~\ref{diffcalculus}; hence
$(\phi_{\mathfrak X},\phi_\Omega)$ is an isomorphism of differential
calculi.
\end{proof}

\section{Application:
 Poisson structures}\label{sect:Poisson}
Having established differential calculi on polyvector fields and mixed differential forms, we now incorporate a Poisson bivector into this framework. We introduce the bigraded Poisson cochain complex and the mixed Poisson chain complex, and show that their cohomology and homology inherit a differential calculus structure.

\subsection{Poisson structure, Poisson homology and cohomology}

This subsection introduces the Poisson cochain complex and the mixed
Poisson chain complex associated with a Poisson algebra. Their
cohomology and homology provide the two spaces underlying the
differential calculus constructed in the next subsection.

\begin{definition}[Poisson structure]
Let $A$ be a bigraded commutative algebra and let
$\epsilon=(r,t)\in\mathbb Z\times\mathbb Z$.
A \textit{Poisson bracket} on $A$ is a bilinear operation
$\{-,-\}:A\otimes A\to A$ with parities specified by $r$ and $t$
in the two gradings, satisfying the following properties
for homogeneous $f,g,h\in A$:
\begin{enumerate}
\item[$(1)$] the bigraded skew commutativity:
\[
\{f,g\}=-(-1)^{(|f|+\epsilon)(|g|+\epsilon)}\{g,f\};
\]

\item[$(2)$] the bigraded Jacobi identity:
\begin{align*}
(-1)^{(|f|+\epsilon)(|h|+\epsilon)}\{f,\{g,h\}\}
+(-1)^{(|g|+\epsilon)(|h|+\epsilon)}\{h,\{f,g\}\}\\
+(-1)^{(|g|+\epsilon)(|f|+\epsilon)}\{g,\{h,f\}\}=0;
\end{align*}

\item[$(3)$] the bigraded Leibniz rule:
\[
\{f,gh\}=\{f,g\}h
+(-1)^{|g|(|f|+\epsilon)}g\{f,h\}.
\]
\end{enumerate}
The bracket is called \textit{even} or \textit{odd}
according as $r$ is even or odd.
\end{definition}

Any bivector field $\pi$ on $A$ defines a bilinear operation
$\{-,-\}:A\otimes A\to A$ by
\begin{align*}
    \{f,g\}:=\pi(df,dg),\ \mathrm{for\ any}\ f,g\in A.
\end{align*}
A bivector field $\pi$ is called a \textit{Poisson bivector field}
if this operation is a Poisson bracket.

Under the degree assumption below, we use the suspension convention
$\pi(df,dg)=[[f,\pi],g]$.
The Schouten--Nijenhuis identities give graded skew-commutativity
and the Leibniz rule.
For $A=\mathcal O(\Sigma)$, 
the Jacobi identity is equivalent to $[\pi,\pi]=0$.
This gives the following proposition.

\begin{proposition}
Let $\pi$ be a bivector field on $A=\mathcal O(\Sigma)$ whose
homogeneous components have even parity degree and even weight.
Then $\pi$ is a Poisson bivector field if and only if $[\pi,\pi]=0$.
\end{proposition}
A bracket $\{-,-\}$ defined in this way is called the
\textit{Poisson bracket} associated with $(A,\pi)$.

\begin{definition}\label{cohomology def}
Suppose that $(A,\pi)$ is a Poisson algebra and
$\deg(\pi)=(0,2)$.
The \textit{Poisson cochain complex} $(\mathrm{CP}^{\bullet,\bullet}
(A),\delta^{\pi}_{\bullet,\bullet})$ is defined as follows. For any  
$(s,t)\in \mathbb{Z}\times\mathbb{Z}$, 
\[
\mathrm{CP}^{s,t}(A):=\mathfrak{X}^{-s,t}(A)
=\{P\in \mathfrak{X}^{\bullet}(A)
\mid\deg(P)=(-s,t)\}.
\]
Thus the upper indices $(s,t)$ on $\mathrm{CP}$ correspond to the actual bidegree $(-s,t)$.  Since the Schouten--Nijenhuis bracket has actual bidegree
$(0,-1)$, the coboundary is
\[
\delta^{\pi}_{s,t}:=[-,\pi]:
\mathrm{CP}^{s,t}(A)\longrightarrow\mathrm{CP}^{s,t+1}(A).
\]
The identity $(\delta^\pi)^2=0$ follows from the shifted Jacobi identity
and $[\pi,\pi]=0$.
The associated cohomology 
$\mathrm{HP}^{s,t}(A)$  
is called \textit{Poisson cohomology}.
\end{definition}

Let $\Sigma$ be a bigraded space and let
$A=(\mathcal{O}(\Sigma),\pi)$ be a Poisson algebra.
Define a linear map $\partial_{\pi}:=[\tilde\iota_{\pi}, \tilde 
d]:\tilde{\Omega}^\bullet(A)\to \tilde{\Omega}^{\bullet}(A)$.
\begin{proposition}\label{boundarypro}
Let \(A=(\mathcal{O}(\Sigma), \pi)\) be a Poisson algebra with $\mathrm{deg}(\pi)=(0,2)$. The linear map \(\partial_\pi\) 
commutes (in the bigraded sense) with \(\tilde{d}\) ,\(\tilde\iota_{\pi}\)
and satisfies \(\partial_\pi \circ \partial_\pi = 0\).
\end{proposition}

\begin{proof} 
By mixed Cartan's formula \eqref{mixedcartanformula},
\[
[\partial_\pi,\tilde\iota_{\pi}] =\tilde\iota_{[\pi, \pi]} = 0,
\]
where the last equality follows from the fact that $\pi$ is a Poisson structure.

For the commutator with $\tilde{d}$, we have
\begin{align*}
    [\partial_\pi, \tilde d] &= \partial_\pi \circ\tilde d 
    +\tilde d \circ \partial_\pi \\ 
    &=\tilde\iota_{\pi} \circ\tilde d\circ\tilde d -\tilde d \circ 
    \tilde\iota_{\pi}\circ\tilde d+\tilde d\circ \tilde\iota_{\pi} 
    \circ\tilde d -\tilde d\circ\tilde d \circ\tilde \iota_{\pi}\\ 
    &= 0.
\end{align*}  
Finally, since the degree of $\pi$ is $(0,2)$, we have
\begin{align*}
  2 \partial_\pi \circ \partial_\pi=[\partial_\pi,\partial_\pi]=[\partial_\pi,
  [\tilde\iota_\pi,\tilde d]]=[\tilde d, [\partial_\pi,\tilde\iota_\pi]]=0,
\end{align*}
which completes the proof.
\end{proof}

Since $\partial_{\pi}$ is a differential, it defines a chain complex for $A$.

\begin{definition}
Let $\Sigma$ be a bigraded space.
Suppose that $A=(\mathcal{O}(\Sigma),\pi)$ is a Poisson algebra, 
and $\mathrm{deg }(\mathrm{\pi})=(0,2)$.
The chain complex of $A$, denoted by 
$(\widetilde{\mathrm{CP}}_{\bullet,\bullet}
(A),\partial_\pi^{\bullet,\bullet})$, is called the {\it  (mixed) Poisson 
chain complex of $A$}, where 
$\widetilde{\mathrm{CP}}_{s,t}(A):=\tilde{\Omega}^{s,-t}(A),$
and the boundary map is given by 
\begin{align*}
\partial^{s,t}_{\pi}:=[\tilde\iota_\pi,\tilde d]:
\widetilde{\mathrm{CP}}_{s,t}(A)\to 
\widetilde{\mathrm{CP}}_{s,t-1}(A).
\end{align*}
The associated homology, denoted by
$\widetilde{\mathrm{HP}}_{s,t}(A)$ 
is called the (\textit{mixed}) \textit{Poisson homology} of $A$.
\end{definition}

\subsection{Differential calculus on Poisson (co)homology}

The operations constructed in Section~\ref{sect:diffcalcul} are compatible with the Poisson coboundary and boundary operators and therefore descend to Poisson cohomology and mixed Poisson homology. We show that the induced wedge product, Schouten--Nijenhuis bracket, contraction, and differential satisfy the differential calculus identities, yielding Theorem~\ref{DC-on-Poisson}.

By the biderivation identity and the shifted Jacobi identity,
the wedge product and the bracket send cocycles to cocycles and send
a product or bracket involving a coboundary and a cocycle to a
coboundary.  Both operations therefore descend to
$\mathrm{HP}^{\bullet,\bullet}(A)$.

\begin{theorem}\label{DC-on-Poisson}
Suppose that $A=(\mathcal{O}(\Sigma),\pi)$ is a Poisson algebra and
$\deg(\pi)=(0,2)$. 
Then
$$\big(\mathrm{HP}^{\bullet,\bullet}(A), 
\widetilde{\mathrm{HP}}_{\bullet,\bullet}(A), 
\wedge,\tilde\iota, [-,-], \tilde d\big)
$$
is a differential calculus.
\end{theorem}

\begin{proof}
In Definition~\ref{diffcalculus}, set
$\mathrm M^{i,j}(A):=\widetilde{\mathrm{HP}}_{-i,-j}(A)$.

(1) 
By the preceding discussion and Theorem~\ref{DC-on-bigraded}, $\mathrm{HP}^{\bullet,\bullet}(A)$ is also a Gerstenhaber algebra.

(2) The map $\tilde{\iota}$: $\mathrm{HP}^{s,t}(A)\otimes
\widetilde{\mathrm{HP}}_{u,v}(A)
\to \widetilde{\mathrm{HP}}_{u-s,v-t}(A)$ is also induced from 
\eqref{mixedcontraction}.
We now check that $\tilde{\iota}$ is well defined on 
$\mathrm{HP}^{\bullet,\bullet}(A)$ and 
$\widetilde{\mathrm{HP}}_{\bullet,\bullet}(A)$.
For any homogeneous $P\in\mathrm{CP}^{s,t}(A)$, the mixed Cartan identity
gives
\[
[\partial_{\pi},\tilde{\iota}_{P}]
=(-1)^{w(P)}\tilde{\iota}_{[P,\pi]}.
\]
If $P$ is a Poisson cocycle, $\tilde\iota_P$ is therefore a graded chain map, so
it sends mixed Poisson cycles to cycles and boundaries to boundaries.  If
$P=\delta^\pi R=[R,\pi]$ is a Poisson coboundary, then
\[
\tilde\iota_P=(-1)^{w(R)}[\partial_\pi,\tilde\iota_R],
\]
so $\tilde\iota_P$ is graded chain null-homotopic and induces the zero map on
mixed Poisson homology.  This proves independence of both the cochain and
chain representatives and yields the asserted action of Poisson
cohomology on mixed Poisson homology.

(3) Finally, Proposition~\ref{boundarypro} gives
$[\tilde d,\partial_\pi]=0$.  Hence $\tilde d$ sends cycles to cycles and
boundaries to boundaries and is well defined on mixed Poisson homology.
Moreover, in $\widetilde{\mathrm{HP}}_{\bullet,\bullet}(A)$, 
if we set $L_P:=[\tilde\iota_P,\tilde d]
=\tilde\iota_P\circ \tilde d-(-1)^{|P||\tilde{d}|}
\tilde d\circ\tilde\iota_P$, 
then by \eqref{mixedcartanformula}, we have:
\[
\tilde\iota_{[P,Q]}=[L_P,\tilde\iota_Q] :={L}_{P}\circ
\tilde\iota_Q-(-1)^{|L_P||Q|}\tilde\iota_Q\circ {L}_P.
\qedhere
\]
\end{proof}

\subsection{Volume
form, unimodularity and the 
Batalin-Vilkovisky algebra}\label{sect:bv}

This subsection recalls bigraded differential calculus and the role of volume forms. We then discuss differential calculus with duality and the resulting Batalin--Vilkovisky structures.

Suppose that $(A,\pi)$ is a Poisson algebra with
$A=\mathcal O(\Sigma)$, and
$\gamma\in\widetilde{\mathrm{CP}}_{u,v}(A)$ for some $u,v\in\mathbb Z$.
We call $\gamma$ a \textit{mixed volume form} if
$\gamma\in\widetilde{\Omega}^{\dim\Sigma}(A)$ and $\mathrm{CP}^{i,j}
(A)\stackrel{\tilde\iota_{(-)}\gamma}{\longrightarrow}
\widetilde{\mathrm{CP}}_{u-i,v-j}(A)$
is an isomorphism of vector spaces.

If $A$ is a Poisson algebra, consider the diagram
\begin{equation}\label{diag:unimodularPoisson}
\begin{array}{ccc}
\mathrm{CP}^{\bullet,\bullet}(A)
&
\xrightarrow{\ \tilde{\iota}_{(-)}\gamma\ }
&
\widetilde{\mathrm{CP}}_{u-\bullet,v-\bullet}(A)
\\
{\scriptstyle [\pi,-]}\,\Big\uparrow
&&
\Big\uparrow\,{\scriptstyle \partial_\pi}
\\
\mathrm{CP}^{\bullet,\bullet-1}(A)
&
\xrightarrow{\ \tilde{\iota}_{(-)}\gamma\ }
&
\widetilde{\mathrm{CP}}_{u-\bullet,v-\bullet+1}(A)
\end{array}
\end{equation}
which need not commute; equivalently, $\gamma$ need not be a Poisson
cycle, that is, $\partial_\pi(\gamma)$ need not vanish. We say that $A$
is {\it unimodular} if there exists a mixed volume form $\gamma$ such
that \eqref{diag:unimodularPoisson} commutes, or equivalently, such that
$\gamma$ is a Poisson cycle. Indeed, if
$\operatorname{div}_{\gamma}(\pi)$ is defined by
$\widetilde{\iota}_{\operatorname{div}_{\gamma}(\pi)}\gamma
=
\widetilde d\bigl(\widetilde{\iota}_{\pi}\gamma\bigr)$,
then $\widetilde d\gamma=0$ and
$\partial_\pi=[\widetilde{\iota}_{\pi},\widetilde d]$ give
$\partial_\pi(\gamma)
=
-\widetilde{\iota}_{\operatorname{div}_{\gamma}(\pi)}\gamma$,
and the fact that contraction with $\gamma$ is an isomorphism yields
\[
\partial_\pi(\gamma)=0
\quad\Longleftrightarrow\quad
\operatorname{div}_{\gamma}(\pi)=0.
\]

We now introduce Lambre's
{\it differential calculus with duality}
(see 
\cite{lambre2010dualite}).
As a consequence, the Gerstenhaber algebra in such a differential calculus 
admits a Batalin--Vilkovisky algebra structure.

\begin{definition}[Lambre]
\label{diff calculus}
A differential calculus $(\mathrm{H}^{\bullet,\bullet},
\mathrm{H}_{\bullet,\bullet}, \cup, \iota, [-,-], d)$ is
called a \textit{differential calculus with duality} if 
there exist integers $u,v$ and an element $\gamma\in \mathrm{H}_{u,v}$
such that
\begin{enumerate}
\item[(a)]  $\iota_{1} \gamma=\gamma$, 
where $1\in \mathrm{H}^{0,0}$ is the unit, $d(\gamma)=0$, and
\item[(b)]  for any $i,j\in \mathbb Z$,
\begin{equation}\label{PD}
\mathrm{PD}(-):=\iota_{(-)}\gamma : \mathrm{H}^{i,j}
\to \mathrm{H}_{u-i,v-j}
\end{equation}
is an isomorphism.
\end{enumerate}
Such an isomorphism $\mathrm{PD}$ is called the 
\textit{Van den Bergh duality} (also called \textit{the
noncommutative Poincar\'e duality}).
\end{definition}

\begin{definition}[Batalin-Vilkovisky algebra]
Suppose that $A$
is a bigraded commutative algebra. A {\it Batalin-Vilkovisky
algebra} structure on $A$ is $(A , \Delta)$
such that
\begin{enumerate}
\item[$(1)$] $\Delta: A^{i,j}\to A^{i,j-1}$ is 
a differential of degree $(0,-1)$, 
that is, $\Delta^2=0$; 
\item[$(2)$] $\Delta$ is a second-order operator,
that is,
\begin{align*}
\Delta(a  b c)
=& \Delta(a b)c+(-1)^{|a||\Delta|}a 
\Delta(b  c)+(-1)^{|b|(|a|+|\Delta|)}b  
\Delta(a   c)\nonumber \\
&- (\Delta a)  b   c-(-1)^{|a||\Delta|}a 
(\Delta b) c-(-1)^{|\Delta|(|a|+|b|)}a  b(\Delta c).
\end{align*}
\end{enumerate}
\end{definition}
The following proposition generalizes the result 
of Getzler in \cite[Proposition 1.2]{getzler1994batalin}.

\begin{proposition}\label{Prop:BVoperator}
A bigraded Batalin-Vilkovisky algebra $A$ is a 
Gerstenhaber algebra equipped with an operator 
$\Delta:A^{\bullet,\bullet}\to A^{\bullet,\bullet-1}$ 
such that $\Delta^2=0$
and that
\begin{align}\label{BVoperator}
[a,b]=(-1)^{|\Delta||a|}\Delta(a b)-(-1)^{|\Delta||a|}\Delta(a) b-a \Delta(b),
\end{align}
for any homogeneous elements $a,b\in A$.
Moreover, in a bigraded Batalin--Vilkovisky algebra, 
$\Delta$ satisfies the identity
$$\Delta[a,b]=[\Delta a,b]+(-1)^{(|a|+|\Delta|)|\Delta|}[a,\Delta b].$$
\end{proposition}

\begin{proof}
First, we prove that $[-,-]$ defined by \eqref{BVoperator} is skew-symmetric:
\begin{align*}
    [a,b]&=(-1)^{|a||\Delta|}
    \Delta(a b)-(-1)^{|\Delta||a|}\Delta(a) b-a \Delta(b)\\
    &=(-1)^{|a||\Delta|+|a||b|}\Delta(ba)-(-1)^{|\Delta||a|+|b|
    (|a|+|\Delta|)}b\Delta(a) - (-1)^{|a|(|\Delta|+| b|)}\Delta(b)a\\
    &=-(-1)^{(|a|+|\Delta|)(|b|+|\Delta|)}[b,a]
\end{align*}
for any homogeneous elements $a,b\in A$.

Second, we prove the compatibility of $\Delta$ with $[-,-]$:
\begin{align*}
&[\Delta a,b]+(-1)^{(|a|+|\Delta|)|\Delta|}[a,\Delta b]\\
=&(-1)^{(|a|+|\Delta|)|\Delta|}\Delta((\Delta a)b)
-(\Delta a)(\Delta b)-\Delta(a(\Delta b))+(\Delta a)(\Delta b)\\
=&(-1)^{(|a|+|\Delta|)|\Delta|}\Delta((\Delta a)b)
-\Delta(a(\Delta b))\\
=&\Delta[a,b],
\end{align*}
since $\Delta$ is a differential.  
  
Next we verify $[-,-]$ is a biderivation of bidegree $(0,-1)$ 
with respect to the product. In fact, 
\begin{align*}
&[a,bc]-[a,b]c-(-1)^{(|a|+|\Delta|)|b|}b[a,c]\\
=&(-1)^{|a||\Delta|}\Delta(abc)-(-1)^{|a||\Delta|}
(\Delta a)bc-a\Delta (bc)\\
&-(-1)^{|a||\Delta|}\Delta(ab)c+(-1)^{|a||\Delta|}
(\Delta a)bc+a\Delta (b)c\\
&-(-1)^{(|a|+|\Delta|)|b|}((-1)^{|a||\Delta|}b\Delta(ac)
-(-1)^{|a||\Delta|}b(\Delta a)c-ba\Delta (c))\\
=&(-1)^{|a||\Delta|}\Delta(abc)-(-1)^{|a||\Delta|}\Delta(abc)\\
=&0
\end{align*}
for any homogeneous elements $a,b,c\in A$.

Finally, we check the Jacobi identity for $[-,-]$.
In fact,
\begin{align*}
    &[[a,b],c]+(-1)^{(|a|+|\Delta|)(|b|+|\Delta|)}[b,[a,c]]\\
    =&(-1)^{(|a|+|b|+|\Delta|)|\Delta|}\Delta([a,b]c)-
    (-1)^{(|a|+|b|+|\Delta|)|\Delta|}\Delta([a,b])c-[a,b]\Delta(c)\\
    &+(-1)^{(|a|+|\Delta|)(|b|+|\Delta|)}((-1)^{|b||\Delta|}\Delta(b[a,c])-
    (-1)^{|b||\Delta|}\Delta(b)[a,c]-b\Delta([a,c]))\\
    =&(-1)^{(|a|+|b|+|\Delta|)|\Delta|}\Delta([a,bc])-
    (-1)^{(|a|+|b|+|\Delta|)|\Delta|}[\Delta a,b]c-(-1)^{(|a|+|\Delta|)
    (|b|+|\Delta|)}b[\Delta a,c])\\
    &-(-1)^{|b||\Delta|}[a,\Delta b]c-(-1)^{|b||\Delta|+(|a|+|\Delta|)
    (|b|+|\Delta|)}\Delta(b)[a,c]\\
    &-[a,b]\Delta(c)-(-1)^{(|a|+|\Delta|)(|b|+|\Delta|)+|\Delta|
    (|\Delta|+|a|)}b[a,\Delta c])\\
    =&[a,[b,c]].
\end{align*}
This proves the claim.
\end{proof}

$\Delta$ is also called the Batalin-Vilkovisky operator, 
or the generator (of the 
Gerstenhaber bracket).
If there exists a degree $(0,-1)$ operator $\Delta$ with 
$\Delta^2=0$ on a Gerstenhaber algebra $(A,[-,-])$ such that 
\eqref{BVoperator} holds, then $A$ is a BV algebra.

Now let $(\mathrm{H}^{\bullet,\bullet},
\mathrm{H}_{\bullet,\bullet},\cup,\iota,[-,-],d,\gamma)$
be a differential calculus with duality, and define the linear operator
$\Delta: \mathrm{H}^{\bullet,\bullet}\to\mathrm{H}^{\bullet,\bullet-1}$ by
$\Delta(a)=-\mathrm{PD}^{-1}\circ d\circ\mathrm{PD}(a)$, for $a\in 
\mathrm{H}^{\bullet,\bullet}$.
Lemma \ref{BVlemma} and Theorem \ref{Thm_Lambre} below are 
generalizations of the work in Lambre's paper \cite{lambre2010dualite}.

\begin{lemma}\label{BVlemma}
Let $(\mathrm{H}^{\bullet,\bullet},
\mathrm{H}_{\bullet,\bullet},\cup,\iota,[-,-],d,\gamma)$
be a differential calculus with duality.
Let $\mathrm{PD}^{-1}=(\iota_{(-)}\gamma)^{-1}$ be the inverse of the duality 
isomorphism.
For any homogeneous elements $P,Q\in\mathrm{H}^{\bullet,\bullet}$ and $\omega\in 
\mathrm{H}_{\bullet,\bullet}$, the following relation holds:
\begin{align}\label{2}
[P,Q]\cup \mathrm{PD}^{-1}(\omega)=&-P\cup \Delta(Q\cup 
\mathrm{PD}^{-1}
(\omega))+(-1)^{w(P)}\Delta(P\cup Q\cup \mathrm{PD}^{-1}(\omega))\nonumber\\ \nonumber
&+(-1)^{p(P)p(Q)+w(Q)(w(P)+1)}Q\cup P\cup \Delta(\mathrm{PD}^{-1}
(\omega))\\
&-(-1)^{p(P)p(Q)+w(Q)w(P)+w(P)+w(Q)}Q\cup \Delta(P\cup \mathrm{PD}^{-1}(\omega)).
\end{align}
\end{lemma}

\begin{proof}
First, we prove that 
$$\mathrm{PD}^{-1}(\iota_P\omega)
    =P\cup \mathrm{PD}^{-1}(\omega).$$
    Since $\mathrm{PD}^{-1}$ is an isomorphism, there exists $R\in 
    \mathrm{H}^{\bullet,\bullet}$, such that $\iota_R\gamma=\omega$, i.e., 
    $\mathrm{PD}^{-1}(\omega)=R$.
    Then  $$\mathrm{PD}^{-1}(\iota_P\omega)
    =\mathrm{PD}^{-1}(\iota_{P\cup R}\gamma)=P\cup R=P\cup \mathrm{PD}^{-1}(\omega).$$

By the above equation, we have $\mathrm{PD}^{-1}(\iota_{[P,Q]}\omega)
=[P,Q]\cup 
\mathrm{PD}^{-1}(\omega).$
Furthermore, by the relation \eqref{calculusformula}, 
$\mathrm{PD}^{-1}
(\iota_{[P,Q]}\omega)$ is written as $\mathrm{PD}^{-1}
([[\iota_P,d],\iota_Q]\omega)$.
Hence, \begin{align*}
&\mathrm{PD}^{-1}(\iota_{[P,Q]}\omega)\\
=&\mathrm{PD}^{-1}\circ\iota_P\circ d\circ\iota_Q(\omega)-
(-1)^{w(P)}\mathrm{PD}^{-1}\circ d\circ\iota_P\circ\iota_Q(\omega)\\
&-(-1)^{p(P)p(Q)+w(Q)(w(P)+1)}
\mathrm{PD}^{-1}\circ\iota_Q\circ\iota_P\circ d(\omega)\\
&+(-1)^{p(P)p(Q)+w(Q)w(P)+w(P)+w(Q)}\mathrm{PD}^{-1}\circ\iota_Q\circ 
d\circ\iota_P(\omega)\\
=&-P\cup \Delta(Q\cup 
\mathrm{PD}^{-1}(\omega))+(-1)^{w(P)}\Delta(P\cup Q\cup \mathrm{PD}^{-1}(\omega))\\
&+(-1)^{p(P)p(Q)+w(Q)w(P)+w(Q)}Q\cup P\cup 
\Delta(\mathrm{PD}^{-1}(\omega))\\
&-(-1)^{p(P)p(Q)+w(Q)w(P)+w(P)+w(Q)}Q\cup \Delta(P\cup 
\mathrm{PD}^{-1}(\omega)).  \end{align*}
This calculation completes the claim.
\end{proof}

\begin{theorem}[Lambre]\label{Thm_Lambre}
Let $(\mathrm{H}^{\bullet,\bullet},\mathrm{H}_{\bullet,\bullet},
\cup,\iota,[-,-],d,\gamma)$
be a differential calculus with duality.
Then $(\mathrm{H}^{\bullet,\bullet},\cup,\Delta)$ is a 
Batalin-Vilkovisky algebra.
\end{theorem}
\begin{proof}
Apply \eqref{2} with $\omega=\gamma$. Since
$\mathrm{PD}^{-1}(\gamma)=1$ and $\Delta(1)=0$, we obtain \begin{equation*}
[P,Q]=(-1)^{w(P)}\Delta(P\cup Q)-(-1)^{w(P)}\Delta(P)\cup Q-P\cup\Delta(Q),
\end{equation*}
which shows that $\Delta$ generates the Gerstenhaber bracket on
$\mathrm{H}^{\bullet,\bullet}$. Moreover, $d^2=0$ gives $\Delta^2=0$.
For $m_P(R):=P\cup R$, \eqref{calculusformula} gives
$
[[\Delta,m_P],m_Q]=(-1)^{w(P)}m_{[P,Q]},
$ and $
[[[\Delta,m_P],m_Q],m_R]=0.
$
Together with $\Delta(1)=0$, this is the second-order condition for $\Delta$.
This completes the proof.
\end{proof}

We next prove that for a bigraded Poisson algebra $(\mathcal{O}(\Sigma),\pi )$,
$(\mathrm {HP}^{\bullet,\bullet}(\mathcal{O}(\Sigma)), 
\widetilde{\mathrm{HP}}_{\bullet,\bullet}(\mathcal{O}(\Sigma)))$ with a unimodular Poisson
structure
forms a differential calculus with duality.

\begin{theorem}\label{thm:UnimodtoBV}
Let $\Sigma$ be a finite-dimensional bigraded vector space, let
$A=\mathcal O(\Sigma)$, and let
$\pi\in\mathfrak X^2(A)$ be a homogeneous Poisson bivector of
bidegree $(0,2)$.  Suppose that
$\gamma\in\widetilde{\mathrm{CP}}_{u,v}(A)$ is a mixed volume form
satisfying
$\widetilde d\gamma=0$
and that $\pi$ is unimodular with respect to $\gamma$.  Then
\[
\bigl(
\mathrm{HP}^{\bullet,\bullet}(A),
\widetilde{\mathrm{HP}}_{\bullet,\bullet}(A),
\cup,\widetilde\iota,[-,-],\widetilde d,[\gamma]
\bigr)
\]
is a differential calculus with duality.  Consequently,
$\mathrm{HP}^{\bullet,\bullet}(A)$ is a Batalin--Vilkovisky algebra.
\end{theorem}

\begin{proof}
By Theorem~\ref{DC-on-Poisson},
$\bigl(
\mathrm{HP}^{\bullet,\bullet}(A),
\widetilde{\mathrm{HP}}_{\bullet,\bullet}(A),
\cup,\widetilde\iota,[-,-],\widetilde d
\bigr)$
is a differential calculus.  It remains to verify that the class
$[\gamma]$ provides the required duality.

Since $\gamma$ is a mixed volume form, contraction with $\gamma$
defines, for every $i,j\in\mathbb Z$, an isomorphism
\[
\operatorname{PD}_\gamma
:=
\widetilde\iota_{(-)}\gamma:
\mathrm{CP}^{i,j}(A)
\xrightarrow{\;\sim\;}
\widetilde{\mathrm{CP}}_{u-i,v-j}(A).
\]
This is part of the definition of a mixed volume form and does not
require a choice of coordinates or a particular expression for
$\gamma$.

The unimodularity of $\pi$ with respect to $\gamma$ means that the
diagram~\ref{diag:unimodularPoisson} commutes.  Equivalently, for every homogeneous
$P\in\mathrm{CP}^{i,j}(A)$,
\[
\partial_\pi
\bigl(\widetilde\iota_P\gamma\bigr)
=
(-1)^j\widetilde\iota_{\delta^\pi(P)}\gamma,
\qquad
\delta^\pi(P)=[P,\pi].
\]
Thus
$\partial_\pi\circ\operatorname{PD}_\gamma
=
\operatorname{PD}_\gamma\circ [\pi,-].$
It follows that $\operatorname{PD}_\gamma$ is an
isomorphism between the Poisson cochain complex with differential $[\pi,-]$ and the shifted mixed
Poisson chain complex.  Passing to cohomology and homology yields an
isomorphism
\[
\operatorname{PD}_\gamma:
\mathrm{HP}^{i,j}(A)
\xrightarrow{\;\sim\;}
\widetilde{\mathrm{HP}}_{u-i,v-j}(A),
\qquad
[P]\longmapsto[\widetilde\iota_P\gamma].
\]

Since $\pi$ is unimodular, one obtains
$\partial_\pi\gamma=0$.
Hence $\gamma$ determines a mixed Poisson homology class
$[\gamma]\in\widetilde{\mathrm{HP}}_{u,v}(A)$.
Moreover, the assumed equality $\widetilde d\gamma=0$ gives
$\widetilde d[\gamma]=[\widetilde d\gamma]=0.$
Therefore $[\gamma]$ satisfies condition~\textup{(a)} in the
definition of a differential calculus with duality, while the
isomorphism $\operatorname{PD}_\gamma$ verifies
condition~\textup{(b)}.  Consequently,
\[
\bigl(
\mathrm{HP}^{\bullet,\bullet}(A),
\widetilde{\mathrm{HP}}_{\bullet,\bullet}(A),
\cup,\widetilde\iota,[-,-],\widetilde d,[\gamma]
\bigr)
\]
is a differential calculus with duality.

By Theorem~\ref{Thm_Lambre}, it endows 
$\mathrm{HP}^{\bullet,\bullet}(A)$
with a Batalin--Vilkovisky algebra structure.
\end{proof}

\begin{example}\label{ex:volumeform}
On 
$A=\mathcal{O}(k^{m||n})$
there exists a mixed volume form 
\begin{align}\label{volumeform}
   \gamma= dx_1 \wedge \cdots \wedge dx_m \wedge y_1^* \wedge \cdots \wedge y_n^*.
\end{align}
Similarly, on
$A^!=\mathcal{O}(k^{n\wedge m})$,
\begin{align}\label{volumeform2}\gamma^!=\xi_1^*
\wedge\cdots\wedge \xi_m^*\wedge d\eta_1\wedge \cdots \wedge d\eta_n
\end{align}
is a mixed volume form.
Let $\pi\in\mathfrak X^2(A)$ and
$\pi^!\in\mathfrak X^2(A^!)$ be Poisson bivectors of bidegree
$(0,2)$.  Assume that they are unimodular with respect to the displayed
volume forms; equivalently,
$\partial_\pi\gamma=0$,
$\partial_{\pi^!}\gamma^!=0$.
Hence, in this case, $\mathrm{HP}^{\bullet,\bullet}(A)$ and
$\mathrm{HP}^{\bullet,\bullet}(A^!)$ are Batalin--Vilkovisky algebras.
\end{example}

\section{Koszul duality
of Poisson structures}\label{sect:koszuldual5}

With the preceding constructions in place, we turn to quadratic Poisson structures and the Poisson complexes associated with them.
Inspired by Shoikhet's observation, Chen et al.
\cite{chen2021poisson}
showed that, for Koszul dual Poisson
spaces $\mathbb{R}^{m|0}$ and $\mathbb{R}^{0|m}$, 
the differential calculus
structures on the corresponding
Poisson cohomologies and homologies
are preserved under Koszul duality.
We prove the analogous result for $k^{m||n}$ and $k^{n\wedge m}$.

\subsection{Poisson
structures and Koszul duality}
This subsection studies quadratic Poisson structures on $k^{m||n}$ and $k^{n\wedge m}$. The discussion is inspired by Shoikhet's result \cite{shoikhet2010Koszul} for the case of $k^{m|0}$ and $k^{0|m}$.
We show that $\pi$ is Poisson if and only if so is its Koszul dual $\pi^!$.

\begin{definition}\label{quadratic Poisson polynomial algebra}
Let $\Sigma$ be a bigraded  space, and let $\{z_1,\cdots z_l\}$ 
be the coordinates of $\Sigma$. A Poisson structure $\pi $ on this space is 
called \textit{quadratic} if it has the form 
\[
\pi=\sum_{i_1,i_2,j_1,j_2} C^{i_1,i_2}_{j_1,j_2} z_{i_1}z_{i_2}
{\partial_{{j_1}}} \wedge {\partial_{{j_2}}} ,
\]
where $C^{i_1,i_2}_{j_1,j_2}\in k$ for any $i_1,i_2,j_1,j_2= 1,\cdots,l.$
It is called \textit{admissible quadratic} if, in addition,
$\deg(\pi)=(0,2)$.
\end{definition}

To simplify the notation,
for $$\mathcal{O}(
k^{m||n})
=\mathrm{Sym}\{x_1,\cdots, x_m, y_1,
\cdots, y_n\}\quad\mbox{and}\quad
\mathcal{O}(
k^{n\wedge m})
=\mathrm{Sym}\{\xi_1,\cdots, \xi_m, 
\eta_1,
\cdots, \eta_n\},$$
we denote their 
generators by $z_i$ and
$\zeta_i$, respectively.
We give a correspondence
between the sets of polyvectors on $\mathcal{O}(
k^{m||n})$ and $\mathcal{O}(
k^{n\wedge m})$,
which exchanges the variables as follows: 
\begin{align}\label{correspondence1}
z_i\leftrightarrow 
{\partial_{\zeta_i}},\quad{\partial_{z_i}}\leftrightarrow
(-1)^{p(z_i)+1}\zeta_i,\quad i=1,\cdots,m+n.
\end{align}

Under the correspondence \eqref{correspondence1},
a quadratic bivector field $\pi=\sum_{i_1,i_2,j_1,j_2}C^{i_1,i_2}_{j_1,j_2}z_{i_1}z_{i_2}{\partial_{{j_1}}} 
\wedge 
{\partial_{{j_2}}}$ on $\mathcal{O}(k^{m||n})$ corresponds to a quadratic bivector field $\pi^!$ on $\mathcal{O}(
k^{n\wedge m})$, and vice versa.

\begin{proposition}\label{quadratic poisson Koszul}
Under the above correspondence, a quadratic bivector
$\pi$ on $k^{m||n}$ whose homogeneous components have even parity
degree is Poisson if and only if the corresponding
bivector $\pi ^{!}$ on $k^{n\wedge m}$ is Poisson.
\end{proposition}

\begin{proof}
By definition, $\pi^!=\phi_{\mathfrak X}(\pi)$.
Both bivectors have weight $2$, and $\phi_{\mathfrak X}$
preserves their parity degrees.
The bracket compatibility \eqref{eq:phi-bracket} gives
\[
[\pi^!,\pi^!]
=[\phi_{\mathfrak X}(\pi),\phi_{\mathfrak X}(\pi)]
=\phi_{\mathfrak X}([\pi,\pi]).
\]
Since $\phi_{\mathfrak X}$ is injective, the right-hand side vanishes if
and only if $[\pi,\pi]=0$.
\end{proof}

\begin{example} \label{example: quadratic Poisson}
Take $m=1,n=3$, consider the dual bigraded  vector space 
$\mathbb{R}^{1||3}$ with coordinates $\{x,y_1,y_2,y_3\}$, endowed with dual quadratic 
Poisson bivector field $\pi=(x^2+y_1y_2)\partial_{y_3}\wedge \partial_{y_3}$.

Following the same values of $m,n$, consider the bigraded  vector space $\mathbb{R}^{3\wedge1}$ 
with coordinates $\{\xi, \eta_1,\eta_2,\eta_3\} $, endowed with a quadratic Poisson
bivector field $\pi^!= \eta_3^2 \partial_{\xi}\wedge \partial_{\xi}+ \eta_3^2 
\partial_{\eta_1}\wedge\partial_{\eta_2}$.
\end{example}

Suppose $(\mathcal{O}(k^{m||n}),\pi)$ is a Poisson algebra.
An admissible quadratic bivector field on $k^{m||n}$ is precisely a
$k$-linear combination of the following types:
\begin{align}\label{piform}
    x_ix_j\partial_{x_a}\wedge\partial_{x_b},\quad y_uy_v\partial_{y_s}\wedge\partial_{y_t},\quad x_iy_s\partial_{x_a}\wedge\partial_{y_t},
\end{align}
where $i,j,a,b=1,\cdots,m$ and $s,t,u,v=1,\cdots, n$.

Suppose that $(\mathcal{O}(k^{n\wedge m}),\pi^!)$ is the 
Poisson algebra.  An admissible quadratic bivector field on
$k^{n\wedge m}$ is precisely a $k$-linear combination of the following
types:
\begin{align}\label{piform2} \xi_i\xi_j\partial_{\xi_a}\wedge\partial_{\xi_b},\quad
\eta_u\eta_v\partial_{\eta_s}
\wedge\partial_{\eta_t},\quad\xi_i\eta_s\partial_{\xi_a}
\wedge\partial_{\eta_t},
\end{align}
where $i,j,a,b=1,\cdots,m$ and $s,t,u,v=1,\cdots, n$.

Thus $\deg(\pi)=\deg(\pi^!)=(0,2)$, and the coboundary maps are
\[
\delta^\pi_{s,t}:\mathrm{CP}^{s,t}(A)\longrightarrow
\mathrm{CP}^{s,t+1}(A),
\qquad
\delta^{\pi^!}_{s,t}:\mathrm{CP}^{s,t}(A^!)\longrightarrow
\mathrm{CP}^{s,t+1}(A^!).
\]

\subsection{Isomorphism of differential calculus}

This subsection compares the differential calculus structures associated with the Koszul-dual admissible quadratic Poisson algebras $A=(\mathcal{O}(k^{m||n}),\pi)$ and $A^!=(\mathcal{O}(k^{n\wedge m}),\pi^!)$
and shows 
that
the two differential calculi are isomorphic, completing the proof of Theorem~\ref{mainthm2}(2).


We start with the following proposition.

\begin{proposition}
\label{thm:poisson-koszul-duality}
Let
$A=\bigl(\mathcal O(k^{m||n}),\pi\bigr)$,
$A^!=\bigl(\mathcal O(k^{n\wedge m}),\pi^!\bigr)$
be Koszul-dual quadratic Poisson algebras, where
$\deg(\pi)=\deg(\pi^!)=(0,2)$.
Then the correspondence in Table~\ref{table:1} induces
bidegree-preserving isomorphisms
\begin{align*}
\phi_{\mathfrak X,*}:
\mathrm{HP}^{\bullet,\bullet}(A)
\xrightarrow{\;\sim\;}
\mathrm{HP}^{\bullet,\bullet}(A^!),\quad
\phi_{\Omega,*}:
\widetilde{\mathrm{HP}}_{\bullet,\bullet}(A)
\xrightarrow{\;\sim\;}
\widetilde{\mathrm{HP}}_{\bullet,\bullet}(A^!).
\end{align*}
\end{proposition}

\begin{proof}
We prove the cohomological and homological statements separately.

First consider the Poisson cochain complexes.  Their homogeneous components are
\begin{align*}
\mathrm{CP}^{s,t}(A)
={}&
\mathrm{Sym}^{-s,t}
\bigl\{
x_1,\ldots,x_m,
y_1,\ldots,y_n,
\partial_{x_1},\ldots,\partial_{x_m},
\partial_{y_1},\ldots,\partial_{y_n}
\bigr\},\\
\mathrm{CP}^{s,t}(A^!)
={}&
\mathrm{Sym}^{-s,t}
\bigl\{
\xi_1,\ldots,\xi_m,
\eta_1,\ldots,\eta_n,
\partial_{\xi_1},\ldots,\partial_{\xi_m},
\partial_{\eta_1},\ldots,\partial_{\eta_n}
\bigr\}.
\end{align*}
The assignments
\[
x_i\longmapsto\partial_{\xi_i},
\qquad
\partial_{x_i}\longmapsto-\xi_i,
\qquad
y_j\longmapsto\partial_{\eta_j},
\qquad
\partial_{y_j}\longmapsto\eta_j
\]
preserve the bidegrees of the generators.  Since both cochain spaces
are free bigraded commutative algebras on the displayed generators,
these assignments extend uniquely to a bidegree-preserving
isomorphism
\[
\phi_{\mathfrak X}:
\mathrm{CP}^{\bullet,\bullet}(A)
\xrightarrow{\;\sim\;}
\mathrm{CP}^{\bullet,\bullet}(A^!).
\]

By the definition of the Koszul dual Poisson bivector,
$\phi_{\mathfrak X}(\pi)=\pi^!$.
Moreover, Proposition~\ref{prop:fourier-compatibility} shows that
$\phi_{\mathfrak X}$ preserves the Schouten--Nijenhuis bracket.
Therefore, for every homogeneous
$P\in\mathrm{CP}^{\bullet,\bullet}(A)$,
\begin{align*}
\phi_{\mathfrak X}\bigl(\delta^\pi(P)\bigr)
=
\phi_{\mathfrak X}([P,\pi])
=
[\phi_{\mathfrak X}(P),\phi_{\mathfrak X}(\pi)]
=
[\phi_{\mathfrak X}(P),\pi^!]
=
\delta^{\pi^!}\bigl(\phi_{\mathfrak X}(P)\bigr).
\end{align*}
Thus the diagram
\begin{equation}\label{diag:poisson-cochain-complexes}
\begin{array}{ccc}
\mathrm{CP}^{\bullet,\bullet}(A)
&
\xrightarrow{\ \phi_{\mathfrak X}\ }
&
\mathrm{CP}^{\bullet,\bullet}(A^!)
\\[6pt]
{\scriptstyle\delta^\pi}\big\downarrow
&&
\big\downarrow{\scriptstyle\delta^{\pi^!}}
\\[6pt]
\mathrm{CP}^{\bullet,\bullet+1}(A)
&
\xrightarrow{\ \phi_{\mathfrak X}\ }
&
\mathrm{CP}^{\bullet,\bullet+1}(A^!)
\end{array}
\end{equation}
commutes.  Hence $\phi_{\mathfrak X}$ is an isomorphism of Poisson
cochain complexes.  It identifies cocycles with cocycles and
coboundaries with coboundaries, and therefore induces an isomorphism
\[
\phi_{\mathfrak X,*}:
\mathrm{HP}^{\bullet,\bullet}(A)
\xrightarrow{\;\sim\;}
\mathrm{HP}^{\bullet,\bullet}(A^!).
\]

We next consider the mixed Poisson chain complexes.  Their homogeneous components are
\begin{align*}
\widetilde{\mathrm{CP}}_{s,t}(A)
={}&
\mathrm{Sym}^{s,-t}
\bigl\{
x_1,\ldots,x_m,
\partial_{y_1},\ldots,\partial_{y_n},
dx_1,\ldots,dx_m,
y_1^*,\ldots,y_n^*
\bigr\},\\
\widetilde{\mathrm{CP}}_{s,t}(A^!)
={}&
\mathrm{Sym}^{s,-t}
\bigl\{
\partial_{\xi_1},\ldots,\partial_{\xi_m},
\eta_1,\ldots,\eta_n,
\xi_1^*,\ldots,\xi_m^*,
d\eta_1,\ldots,d\eta_n
\bigr\}.
\end{align*}
The assignments
\[
x_i\longmapsto\partial_{\xi_i},
\qquad
\partial_{y_j}\longmapsto\eta_j,
\qquad
dx_i\longmapsto\xi_i^*,
\qquad
y_j^*\longmapsto d\eta_j
\]
extend uniquely to a bidegree-preserving isomorphism
\[
\phi_\Omega:
\widetilde{\mathrm{CP}}_{\bullet,\bullet}(A)
\xrightarrow{\;\sim\;}
\widetilde{\mathrm{CP}}_{\bullet,\bullet}(A^!).
\]
Equivalently, the correspondence identifies the differentials of the
generators by
\[
\widetilde d(x_i)=dx_i
\longmapsto
\xi_i^*
=
\widetilde d(\partial_{\xi_i}),
\quad
\widetilde d(\partial_{y_j})=y_j^*
\longmapsto
d\eta_j
=
\widetilde d(\eta_j).
\]
It follows that
$\phi_\Omega\circ\widetilde d
=
\widetilde d\circ\phi_\Omega.$

The mixed Poisson boundary operators are
$\partial_\pi
=
[\widetilde\iota_\pi,\widetilde d]$,
$\partial_{\pi^!}
=
[\widetilde\iota_{\pi^!},\widetilde d]$.
By Proposition~\ref{prop:fourier-compatibility},
$\phi_\Omega$ is compatible with both the mixed differential and
contraction.  Since $\phi_{\mathfrak X}(\pi)=\pi^!$, for every
$\omega\in\widetilde{\mathrm{CP}}_{\bullet,\bullet}(A)$ one obtains
\begin{align*}
\phi_\Omega\bigl(\partial_\pi\omega\bigr)
=
\phi_\Omega
\bigl([\widetilde\iota_\pi,\widetilde d]\omega\bigr)
=
[\widetilde\iota_{\phi_{\mathfrak X}(\pi)},\widetilde d]\,
\phi_\Omega(\omega)
=
[\widetilde\iota_{\pi^!},\widetilde d]\,
\phi_\Omega(\omega)
=
\partial_{\pi^!}\bigl(\phi_\Omega(\omega)\bigr).
\end{align*}
Here no additional sign occurs because $\phi_{\mathfrak X}$ and
$\phi_\Omega$ preserve the bidegrees and hence preserve the graded
commutator signs.  Thus the diagram
\[
\begin{array}{ccc}
\widetilde{\mathrm{CP}}_{\bullet,\bullet}(A)
&
\xrightarrow{\ \phi_\Omega\ }
&
\widetilde{\mathrm{CP}}_{\bullet,\bullet}(A^!)
\\[6pt]
{\scriptstyle\partial_\pi}\big\downarrow
&&
\big\downarrow{\scriptstyle\partial_{\pi^!}}
\\[6pt]
\widetilde{\mathrm{CP}}_{\bullet,\bullet-1}(A)
&
\xrightarrow{\ \phi_\Omega\ }
&
\widetilde{\mathrm{CP}}_{\bullet,\bullet-1}(A^!)
\end{array}
\]
commutes.  Consequently, $\phi_\Omega$ is an isomorphism of mixed
Poisson chain complexes.  It identifies cycles with cycles and
boundaries with boundaries, and therefore induces an isomorphism
\[
\phi_{\Omega,*}:
\widetilde{\mathrm{HP}}_{\bullet,\bullet}(A)
\xrightarrow{\;\sim\;}
\widetilde{\mathrm{HP}}_{\bullet,\bullet}(A^!).
\]
This proves both assertions.
\end{proof}

We are now ready to show 
Theorem \ref{mainthm2} (2), which is rephrased
as follows for the reader's convenience.

\begin{theorem}
\label{mainthm22}
Let
\(
A=(\mathcal O(k^{m||n}),\pi),\ \text{and}\ 
A^!=(\mathcal O(k^{n\wedge m}),\pi^!)
\)
be Koszul dual admissible quadratic Poisson algebras.
Then
$(\mathrm{HP}^{\bullet,\bullet}(A),
\widetilde{\mathrm{HP}}_{\bullet,\bullet}(A))$ and 
$(\mathrm{HP}^{\bullet,\bullet}( A^!),
\widetilde{\mathrm{HP}}_{\bullet,\bullet}( A^!))$ are differential calculi.
Under the correspondence
given in Table \ref{table:1},
the two differential calculus
structures are isomorphic.
\end{theorem}

\begin{proof}
The existence of the two differential calculi follows from
Theorem~\ref{DC-on-Poisson}. It remains to prove that the induced maps
intertwine their calculus operations.

By Proposition
\ref{thm:poisson-koszul-duality}, we gives an isomorphism 
$\phi_{\mathfrak{X},*}:\mathrm{HP}^{\bullet,\bullet}(A)\cong 
\mathrm{HP}^{\bullet,\bullet}(A^!)$ as cohomologies.
Equations~\eqref{eq:phi-product} and~\eqref{eq:phi-bracket} descend to
the wedge product and Schouten--Nijenhuis bracket on Poisson cohomology.
By Proposition \ref{thm:poisson-koszul-duality}, there is an isomorphism between 
$\widetilde{\mathrm{HP}}_{\bullet,\bullet}(A) $ and 
$ \widetilde{\mathrm{HP}}_{\bullet,\bullet}(A^!)$.

Equations~\eqref{eq:phi-contraction} and
\eqref{eq:phi-differential} descend to the contraction and mixed
differential on homology.  Consequently,
\[\big(\mathrm{HP}^{\bullet,\bullet}(A), 
\widetilde{\mathrm{HP}}_{\bullet,\bullet}(A)\big)\  
\mathrm{and}\  \big(\mathrm{HP}^{\bullet,\bullet}
(A^!), \widetilde{\mathrm{HP}}_{\bullet,\bullet}
(A^!)\big)\]
are isomorphic.
\end{proof}

\subsection{Isomorphism of Batalin–Vilkovisky algebras}

Having established the isomorphism of differential calculi, we now incorporate the duality data needed to obtain Batalin--Vilkovisky structures. For the canonical mixed volume forms introduced in Section~\ref{sect:Poisson}, we prove that the correspondence preserves unimodularity, which
induces an isomorphism of Batalin--Vilkovisky algebras, completing the proof of Theorem~\ref{mainthm3}(2).

\begin{proposition}
\label{thm:thirdtheorem}
Let $A=(\mathcal{O}(k^{m||n}),\pi)$ be an admissible quadratic Poisson
algebra and let
$A^!=(\mathcal O(k^{n\wedge m}),\pi^!)$ be its Koszul dual.
Given the volume forms $\gamma$ \eqref{volumeform} on $k^{m||n}$ and 
\eqref{volumeform2} $\gamma^!$ on $k^{n\wedge m}$,
$A$
is unimodular with respect to $\gamma$ if and only if its Koszul dual $A^!$ is unimodular with respect to $\gamma^!$.
In this case, 
the following diagram
\begin{equation}\label{commutative diagram}
\begin{array}{ccc}
\mathrm{HP}^{\bullet,\bullet}(A)
&
\xrightarrow{\ \cong\ }
&
\widetilde{\mathrm{HP}}_{-n-\bullet,m-\bullet}(A)
\\[6pt]
{\scriptstyle\cong}\big\downarrow
&&
\big\downarrow{\scriptstyle\cong}
\\[6pt]
\mathrm{HP}^{\bullet,\bullet}(A^!)
&
\xrightarrow{\ \cong\ }
&
\widetilde{\mathrm{HP}}_{-n-\bullet,m-\bullet}(A^!)
\end{array}
\end{equation}
commutes.
\end{proposition}

\begin{proof}
First, we show that
$A=(\mathcal{O}(k^{m||n}), \pi)$ is unimodular if 
and only if $A^!=(\mathcal{O}(k^{n\wedge m}),\pi^{!})$ is unimodular.

Recall the Poisson (co)chain complexes
and let
$$
\gamma=dx_1 \wedge \cdots \wedge dx_m \wedge y_1^* \wedge \cdots 
\wedge y_n^*\quad\mbox{and}\quad
\gamma^!=\xi_1^*
\wedge\cdots\wedge \xi_m^*\wedge d\eta_1\wedge \cdots \wedge d\eta_n.
$$
Then under the correspondence of Table \ref{table:1}, we have
\begin{equation}\label{idoftwodifferentialcalculus}
\begin{array}{ccc}
x_i \leftrightarrow \partial_{\xi_i},
&
y_j \leftrightarrow \partial_{\eta_j},
&
dx_i \leftrightarrow \xi_i^*,
\\[4pt]
\partial_{x_i} \leftrightarrow -\xi_i,
&
\partial_{y_j} \leftrightarrow \eta_j,
&
y_j^* \leftrightarrow d\eta_j.
\end{array}
\end{equation}
The chosen normalization sends $\gamma$ to $\gamma^!$, and the diagram
\begin{equation}\label{cd_fourhomology}
\begin{array}{ccc}
\mathrm{CP}^{\bullet,\bullet}(A)
&
\xrightarrow{\ \tilde\iota_{(-)}\gamma\ }
&
\widetilde{\mathrm{CP}}_{-n-\bullet,m-\bullet}(A)
\\[4pt]
{\scriptstyle\cong}\big\downarrow
&&
\big\downarrow{\scriptstyle\cong}
\\[4pt]
\mathrm{CP}^{\bullet,\bullet}(A^!)
&
\xrightarrow{\ \tilde\iota_{(-)}\gamma^!\ }
&
\widetilde{\mathrm{CP}}_{-n-\bullet,m-\bullet}(A^!)
\end{array}
\end{equation}
commutes.
Thus $\gamma$ is a Poisson cycle of $A$ if 
and only if $\gamma^!$ is a Poisson
cycle of $A^!$, which proves the claim.

Next, we show that the following diagram
\begin{equation}\label{diag:four_iso}
\begin{array}{ccc}
\mathrm{HP}^{\bullet,\bullet}(A)
&
\xrightarrow{\ \cong\ }
&
\widetilde{\mathrm{HP}}_{-n-\bullet,m-\bullet}(A)
\\[6pt]
{\scriptstyle\cong}\big\downarrow
&&
\big\downarrow{\scriptstyle\cong}
\\[6pt]
\mathrm{HP}^{\bullet,\bullet}(A^!)
&
\xrightarrow{\ \cong\ }
&
\widetilde{\mathrm{HP}}_{-n-\bullet,m-\bullet}(A^!)
\end{array}
\end{equation}
commutes.
The two vertical isomorphisms are given by 
Theorem~\ref{thm:poisson-koszul-duality},
and the two horizontal isomorphisms are given by 
Theorem~\ref{thm:UnimodtoBV}.
The commutativity of the diagram \eqref{diag:four_iso} 
follows from the chain level commutative diagram \eqref{cd_fourhomology}.
\end{proof}

By Theorem \ref{thm:UnimodtoBV}, the sextuple $(\mathrm{HP}^{\bullet,\bullet}
(A),\widetilde{\mathrm{HP}}_{\bullet,\bullet}(A), \cup, \tilde\iota, 
[-,-],\tilde d)$ with $\gamma$ (and the dual side with $\gamma^!$) forms 
a differential calculus with duality.
In summary, we obtain the following corollary.

\begin{corollary}\label{thm:Batalin-Vilkovisky}
Let $A =(\mathcal{O}(k^{m||n}),\pi)$ 
and $A^!=(\mathcal{O}(k^{n\wedge m}),\pi^!)$ be 
unimodular admissible quadratic Poisson algebras as above.
Then $\mathrm{HP}^{\bullet,\bullet}(A)$ and 
$\mathrm{HP}^{\bullet,\bullet}(A^!)$ are 
Batalin-Vilkovisky algebras.
\end{corollary}

\begin{example}[Example~\ref{ex:volumeform} continued]
    Let $\gamma$ be the mixed volume form in Example~\ref{ex:volumeform} and let
\[
\begin{aligned}
\pi={}&
\sum_{1\leq a<b\leq m}\sum_{i,j=1}^{m}
A_{ab}^{ij}x_i x_j\,
\partial_{x_a}\wedge\partial_{x_b}
+\sum_{a,i=1}^{m}\sum_{s,r=1}^{n}
B_{as}^{ir}x_i y_r\,
\partial_{x_a}\wedge\partial_{y_s} \\
&+\sum_{1\leq s<t\leq n}\sum_{r,q=1}^{n}
D_{st}^{rq}y_r y_q\,
\partial_{y_s}\wedge\partial_{y_t}.
\end{aligned}
\]
We set $A_{ba}^{ij}:=-A_{ab}^{ij}$ for $a<b$ and
$D_{ts}^{rq}:=D_{st}^{rq}$ for $s<t$, with
$A_{aa}^{ij}=D_{ss}^{rq}=0$.
A direct computation of $\operatorname{div}_{\gamma}(\pi)$ shows that $\pi$
is unimodular with respect to $\gamma$ if and only if
\[
\sum_{a=1}^{m}\bigl(A_{ab}^{ai}+A_{ab}^{ia}\bigr)
+\sum_{s=1}^{n}B_{bs}^{is}=0,
\qquad 1\leq b,i\leq m,
\]
and
\[
\sum_{a=1}^{m}B_{at}^{ar}
+\sum_{s=1}^{n}\bigl(D_{st}^{sr}-D_{st}^{rs}\bigr)=0,
\qquad 1\leq t,r\leq n.
\]
Under the correspondence in Table~\ref{table:1}, the same coefficient conditions
characterize the unimodularity of $\pi^!$ with respect to $\gamma^!$.
\end{example}

We are now ready to
prove Theorem \ref{mainthm3} (2).

\begin{theorem}[Theorem~\ref{mainthm3} (2)]
\label{mainthm33}
Let
$A=(\mathcal O(k^{m||n}),\pi)$
 and 
$A^!=(\mathcal O(k^{n\wedge m}),\pi^!)$
be Koszul-dual admissible quadratic Poisson algebras, equipped with the mixed volume
forms $\gamma$ and $\gamma^!$ in \eqref{volumeform} and
\eqref{volumeform2}. Then $\pi$ is unimodular if and only if $\pi^!$ is
unimodular. In this case,
\[
\mathrm{HP}^{\bullet,\bullet}(A)
\cong
\mathrm{HP}^{\bullet,\bullet}(A^!)
\]
is an isomorphism of Batalin--Vilkovisky algebras.
\end{theorem}

\begin{proof}
By Proposition~\ref{thm:thirdtheorem},
the right vertical isomorphism preserves the mixed differential as well as
the mixed volume form,
so the two differential calculi with duality
$$
\big(\mathrm{HP}^{\bullet,\bullet}(A),
\widetilde{\mathrm{HP}}_{\bullet,\bullet}(A)\big)\;
\mbox{and} \;\big (\mathrm{HP}^{\bullet,\bullet}
(A^{!}),\widetilde{\mathrm{HP}}_{\bullet,\bullet}(A^!)\big)
$$
are isomorphic.
The result now follows from Corollary \ref{thm:Batalin-Vilkovisky}.

It remains to verify compatibility with the Batalin--Vilkovisky operators.
For any homogeneous $[P]\in \mathrm{HP}^{\bullet,\bullet}(A)$, the commutativity of
diagram~\ref{diag:four_iso} and the compatibility of $\phi_{\Omega,*}$ with $\widetilde d$ give
\[
\begin{aligned}
\phi_{\mathfrak X,*}\bigl(\Delta_A[P]\bigr)
&=-
  \phi_{\mathfrak X,*}\circ\mathrm{PD}_{\gamma}^{-1}\circ
  \widetilde d\circ\mathrm{PD}_{\gamma}([P]) \\
&=-
  \mathrm{PD}_{\gamma^!}^{-1}\circ\phi_{\Omega,*}
  \circ\widetilde d\circ\mathrm{PD}_{\gamma}([P]) \\
&=-
  \mathrm{PD}_{\gamma^!}^{-1}\circ\widetilde d\circ
  \mathrm{PD}_{\gamma^!}\circ\phi_{\mathfrak X,*}([P]) \\
&=\Delta_{A^!}\bigl(\phi_{\mathfrak X,*}([P])\bigr).
\end{aligned}
\]
Thus $\phi_{\mathfrak X,*}$ is an isomorphism of
Batalin--Vilkovisky algebras.
\end{proof}

\begin{ack}
We thank Xiaojun Chen for his guidance throughout this work.
We also thank Si Li
for assistance with the signs
(see \eqref{signconv})
and for drawing our attention to \cite{deligne1999notes}. We thank Farkhod Eshmatov for inviting us to New Uzbekistan University.
This work was supported by the NSFC (12271377 and 12261131498).
\end{ack}


\end{document}